\documentclass{amsart}

\usepackage{amsthm,amsfonts,amsmath,amssymb,latexsym,epsfig}
\usepackage{upref,amssymb,eucal}
\usepackage[matrix,arrow,curve,cmtip]{xy}

\DeclareMathOperator{\Spec}{Spec}
\DeclareMathOperator{\Proj}{Proj}
\DeclareMathOperator{\Hom}{Hom}
\newcommand{\longlabelmap}[1]{{\,\buildrel #1\over\longrightarrow\,}}
\newcommand{\pr}{{\mathrm{pr}}}

\begin{document}
\title{Differential forms on arithmetic jet spaces}
\author{James Borger and Alexandru Buium}
\def \bW{{\mathbb W}}
\def \ZN{\bZ[1/N,\zeta_N]}
\def \tcS{\tilde{S}}
\def \Z{{\mathbb Z}}
\def \cI{I}
\def \cU{\mathcal U}
\def \cF{\mathcal F}
\def \tcM{\tilde{M}}
\def \cI{I}
\def \tcI{\tilde{I}}
\def \cK{\mathcal K}
\def \cH{\mathcal H}
\def \cD{\mathcal D}
\def \cE{E}
\def \cP{\mathcal P}
\def \cA{A}
\def \cV{\mathcal V}
\def \cM{M}
\def \cS{\mathcal S}
\def \cN{\mathcal N}
\def \tcM{\tilde{M}}
\def \cG{\mathcal G}
\def \cB{\mathcal B}
\def \cJ{\mathcal J}
\def \tG{\tilde{G}}
\def \cF{\mathcal F}
\def \h{\hat{\ }}
\def \hp{\hat{\ }}
\def \tS{\tilde{S}}
\def \tP{\tilde{P}}
\def \tA{\tilde{A}}
\def \tX{\tilde{X}}
\def \tT{\tilde{T}}
\def \tE{\tilde{E}}
\def \tV{\tilde{V}}
\def \tC{\tilde{C}}
\def \tI{\tilde{I}}
\def \tU{\tilde{U}}
\def \tG{\tilde{G}}
\def \tu{\tilde{u}}
\def \tx{\tilde{x}}
\def \tL{\tilde{L}}
\def \tY{\tilde{Y}}
\def \d{\delta}
\def \bZ{{\mathbb Z}}
\def \bV{{\bf V}}
\def \bE{{\bf E}}
\def \bC{{\bf C}}
\def \bO{{\bf O}}
\def \bR{{\bf R}}
\def \bA{{\mathbb A}}
\def \bZg{\bZ_{\geq 0}}
\def \bB{{\bf B}}
\def \cO{\mathcal O}
\def \ra{\rightarrow}
\def \bX{{\bf X}}
\def \bH{{\bf H}}
\def \bS{{\bf S}}
\def \bF{{\mathbb F}}
\def \bN{{\bf N}}
\def \bK{{\bf K}}
\def \bE{{\bf E}}
\def \bF{{\bf F}}
\def \bB{{\bf B}}
\def \bQ{{\bf Q}}
\def \bd{{\bf d}}
\def \bY{{\bf Y}}
\def \bU{{\bf U}}
\def \bL{{\bf L}}
\def \bQ{{\mathbb Q}}
\def \bP{{\bf P}}
\def \bR{{\bf R}}
\def \bC{{\mathbb C}}
\def \bD{{\bf D}}
\def \bM{{\bf M}}
\def \bG{{\bf G}}
\def \bP{{\bf P}}
\def \tM{\tilde{M}}
\def \tS{\tilde{S}}

\newtheorem{THM}{{\!}}[section]
\newtheorem{THMX}{{\!}}
\renewcommand{\theTHMX}{}
\newtheorem{theorem}{Theorem}[section]
\newtheorem{corollary}[theorem]{Corollary}
\newtheorem{lemma}[theorem]{Lemma}
\newtheorem{proposition}[theorem]{Proposition}
\theoremstyle{definition}
\newtheorem{definition}[theorem]{Definition}
\theoremstyle{remark}
\newtheorem{remark}[theorem]{Remark}
\theoremstyle{assumption}
\newtheorem{assumption}[theorem]{\bf Assumption}
\newtheorem{example}[theorem]{\bf Example}
\numberwithin{equation}{section}
\address{The Australian National University\\ Canberra, ACT 0200, Australia}
\email{james.borger@anu.edu.au}
\address{University of New Mexico \\ Albuquerque, NM 87131, USA}
\email{buium@math.unm.edu}
\maketitle

\begin{abstract}
We study derivations and differential forms on the arithmetic jet spaces of smooth schemes, relative to several primes.
As  applications we give a new interpretation of arithmetic Laplacians and we discuss the de Rham cohomology of some specific arithmetic jet spaces.
\end{abstract}

\section{Introduction}

Arithmetic jet spaces with respect to a single prime $p$ (also called $p$-jet spaces) were introduced in
\cite{char} and further studied in a series of papers; see \cite{book} and the bibliography therein. A
multiple-prime generalization of these spaces was introduced independently in \cite{laplace} and
\cite{borger}. As explained in \cite{book,laplace}, these arithmetic jet spaces can be viewed as an
arithmetic analogue of usual jet spaces in differential geometry and classical mechanics; the role of the
derivatives with respect to various directions is played, in the arithmetic setting, by Fermat quotient
operators with respect to various primes. In particular the functions on arithmetic jet spaces with
respect to several primes can be viewed as arithmetic analogues of classical partial differential
equations on manifolds. See \cite{laplace}, for instance, for an arithmetic analogue of Laplacians.
Alternatively, arithmetic jet spaces can be viewed as a realization, within classical algebraic geometry,
of what one might call absolute geometry, or geometry over the field with one element. See
\cite{borger, f1} or the introduction to \cite{book}.

As explained in \cite{book}, pp.\ 88-100, in the case of one prime, the tangent (and cotangent) bundles
of arithmetic jet spaces carry some remarkable structures that are analogous to structures appearing in
classical mechanics. In this paper we extend this to the case of several primes and then, as
applications, we give a new interpretation of the arithmetic Laplacians in \cite{laplace}, and we discuss
the de Rham cohomology of certain specific arithmetic jet spaces.

Indeed in the paper \cite{laplace}, one considers arithmetic jet spaces $Y=\cJ^r_{\cP}(X)$ of order $r
\in \bZ_{\geq 0}^d$ attached to smooth schemes $X$ over $\bZ$, with respect to a finite set of primes
$\cP=\{p_1,\dots ,p_d\}$, and in case $X$ is a one dimensional group scheme, one constructs arithmetic
analogues of Laplacians. These arithmetic Laplacians are constructed as families $(f_1,\dots ,f_d)$ where
each $f_k$ is a formal function on the completion of $Y$ along $Y_{p_k}:=Y \otimes (\bZ/p_k\bZ)$; these
formal functions are required to be ``analytically continued" along the zero section $Z$ of $X$ in the
sense that there exists a formal function $f_0$ on the completion  of $Y$ along $Z$ such that $f_0$ and
$f_k$ coincide on the completion of $Y$ along $Z\cap Y_{p_k}$ for each $k$. In this paper we
want to revisit the idea of analytic continuation between different primes. Indeed we will show that:

\begin{theorem}
The
arithmetic Laplacians $(f_1,\dots ,f_d)$ of \cite{laplace} have the property that the $1$-forms
$df_1,\dots ,df_d$  extend to $1$-forms defined on the whole of $Y$ and these extended $1$-forms all coincide.
\end{theorem}

See Theorems \ref{tmu} and \ref{tme} below, and the discussion preceding them, for the precise statement of the result and the precise definitions of arithmetic Laplacians.
So the
arithmetic Laplacians appear as primitives, existing only on formal neighborhoods of certain divisors, of
a $1$-form that exists on the whole of the arithmetic jet space $Y$. This provides an
alternative view on analytic continuation between different primes and has a global flavor, as opposed
to the formally local flavor of the one in \cite{laplace}.

Here is the plan of the paper. In section 2 we formulate this concept of analytic continuation based on
globally defined differential forms; this makes sense on any scheme $Y$, not necessarily on arithmetic
jet spaces. In section 3 we review and complement the main concepts in \cite{laplace,borger} related to
arithmetic jet spaces. In sections 4 and 5 we extend the discussion in \cite{book} on the tangent and
cotangent bundles of arithmetic jet spaces from the case of one prime to the case of several primes. In
sections 6, 7, and 8 we apply these concepts to the case when $X$ is, respectively, the multiplicative
group, an affine open set of an elliptic curve, or an affine open set of a K3 surface attached, via the
Kummer construction, to a product of elliptic curves. We will perform ``analytic continuation" on the
arithmetic jet spaces of such schemes $X$, between various primes, based on globally defined differential
forms, and as a by-product, we will derive various de~Rham-style consequences for the arithmetic jet
spaces in question. We have chosen to restrict ourselves to the case when $X$
is affine because this simplifies the exposition and also captures the essential points of the theory.
However our theory can be extended to the case when $X$ is not necessarily affine.

\begin{remark}
 One should call the attention upon the fact that our arithmetic jet spaces are objects completely different from Vojta's \cite{vojta}; indeed Vojta's jet spaces are constructed using Hasse-Schmidt derivations (which are morally  ``differentiations in the geometric direction") while our jet spaces are constructed using Fermat quotients (which are morally ``differentiations in the arithmetic directions".)
 \end{remark}

\bigskip

{\bf Acknowledgment} While writing this paper the first author was partially supported by
ARC grant DP0773301, and the second was partially supported by NSF grants DMS-0552314 and DMS-0852591.
We are indebted to D. Bertrand and J-B. Bost for useful conversations on transcendence and formal functions.

\bigskip
\bigskip

\section{Analytic continuation between primes via differential forms}
For any ring $B$ (respectively scheme $Y$) we denote by $\Omega_{B}$ (respectively $\Omega_{Y}$)
the $B$-module (respectively the sheaf on $Y$) of K\"{a}hler differentials of $B$ (respectively $Y$) over $\bZ$. We denote by $$T_{B}=Hom_B(\Omega_{B},B)=Der(B,B)$$
the dual of $\Omega_{B}$; we denote by $T_{Y}$ the dual sheaf of $\Omega_{Y}$, which we refer to as the {\it tangent sheaf} of $Y$. Also for $i \geq 0$ we let $\Omega_B^i=\wedge^i \Omega_B$, $\Omega^i_Y=\wedge^i\Omega_Y$
denote the exterior powers of $\Omega_B$, $\Omega_Y$ respectively. Elements of $H^0(Y,\Omega_Y^i)$ are referred to as $i$-{\it forms} on $Y$. Also we have at our disposal de Rham complexes
$(\Omega_B^*,d)$ and $(\Omega^*_Y,d)$ respectively and we have the usual notions of closed and exact forms.
Recall that if $B$, $Y$ are smooth over $\bZ$ (or more generally over a ring of fractions of $\bZ$) of relative dimension $m$ then $\Omega_B^i$, $\Omega_Y^i$ are locally free of rank $\left( \begin{array}{c} m\\ i \end{array}\right)$; by a  {\it volume form} on $Y$ we will mean an $m$-form on $Y$ that is invertible (i.e. is a basis of
$\Omega_Y^m$). Volume forms are, of course, closed.

If $M$ is a module over a ring $B$ and  if $I$ is an ideal in $B$ we denote by
$$M^{\widehat{I}}=\lim_{\leftarrow} M/I^nM$$
the $I$-adic completion of $M$. We say that $M$ is $I$-{\it adically complete} if the map $M \ra M^{\widehat{I}}$ is an isomorphism.
For any Noetherian scheme $Y$ and any closed subscheme $Z \subset Y$ we denote by $Y^{\widehat{Z}}$ the formal scheme obtained by  completing $Y$ along $Z$. If $I$ is the ideal sheaf of $Z$ in $Y$ we also write
$Y^{\widehat{I}}$ in place of $Y^{\widehat{Z}}$. We let $Z_n \subset Y$ be the closed subscheme with ideal
$I^n$.
We let
$$\Omega^i_{Y^{\widehat{I}}}=\lim_{\leftarrow}
\Omega^i_{Z_n},\ \ \ T_{Y^{\widehat{I}}}=\lim_{\leftarrow}
T_{Z_n}.$$

Elements of $H^0(Y^{\widehat{I}},\Omega^i_{Y^{\widehat{I}}})$ are referred to as $i$-{\it forms} on $Y^{\widehat{I}}$. Of course, $H^0(Y^{\widehat{I}},\Omega^0_{Y^{\widehat{I}}})$ is just the space of formal functions $\cO(Y^{\widehat{I}})$. Again, we have  at our disposal an obvious deRham complex and a notion of closed and exact forms. If $Y$ is smooth over $\bZ$ (or a ring of fractions of $\bZ$) we have a concept  of volume form.

We will use the following basic terminology:

\begin{definition} An   $(i-1)$-form $\eta \in H^0(Y^{\widehat{I}}, \Omega^{i-1}_{Y^{\widehat{I}}})$ is an $I$-{\it adic primitive} of an $i$-form $\nu \in H^0(Y,\Omega_Y^i)$
if $\nu=d \eta$ in $H^0(Y^{\widehat{I}},\Omega^{i}_{Y^{\widehat{I}}})$.
An $i$-form $\nu \in H^0(Y,\Omega_Y^i)$ is called $I$-{\it adically exact} if it has a $I$-adic primitive.
\end{definition}

 {\it In what follows, throughout the paper, we consider a finite set of primes $\cP=\{p_1,...,p_d\} \subset \bZ$ and we will denote by $A_0=S^{-1}\bZ$ a ring of fractions of $\bZ$ with respect to a multiplicative system $S$ of integers coprime to the primes in $\cP$.}

 \medskip

Let $Y$ be a smooth affine scheme of finite type over $A_0$.
We will be interested in the following examples of subschemes $Z \subset Y$. As before we denote by $I$ the ideal defining $Z$.

{\it 1) Vertical case:} $Z$ is defined by the ideal $I=(p_1...p_d)$. In this case
$$Y^{\widehat{p_1...p_d}}=\coprod_{k=1}^d Y^{\widehat{p_k}}$$
and hence
$$H^0(Y^{\widehat{p_1...p_d}},\Omega^{i}_{Y^{\widehat{p_1...p_d}}})=\prod_{k=1}^d
H^0(Y^{\widehat{p_k}},\Omega^{i}_{Y^{\widehat{p_k}}}).$$
In this case we shall use the phrases  $\cP$-{\it adically exact} and  $\cP$-{\it adic primitive} in place of  $I$-adically exact and  $I$-adic primitive respectively.
If $Y$ is smooth over $A_0$ then $\Omega^i_{Y^{\widehat{p_k}}}$ are locally free.

{\it 2) Horizontal case:} $Z$ is the image in $Y$ of an $A_0$-point $P \in Y(A_0)$, i.e. of a section $P:\Spec A_0\ra Y$ of $Y \ra \Spec A_0$ . In this case the ideal $I$ is the kernel of $P^*:\cO(Y)\ra A_0$; by abuse of notation we denote the ideal $I$, again, by $P$. If $\cS \subset Y(A_0)$ is a set of
points we use the phrase {\it $\cS$-addically  exact} to mean $P$-adically exact for all $P \in \cS$.
Let now $P\in Y(A_0)$ be a {\it uniform point}
(in the sense of \cite{laplace}, Definition 2.25); recall that this means that there exists a
Zariski open set $Y_1 \subset Y$ containing the image of $P$ and possessing   an \'{e}tale morphism $Y_1 \ra \bA^N$ to an affine space $\bA^N=\Spec A_0[T]$, where $T=\{T_j\}$ is a tuple of variables.  If this is the case then one can choose  $Y_1$ and  $T$ such that   the ideal in
$\cO(Y_1)$ of the image of $P:\Spec A_0 \ra Y_1$, i.e. the kernel of $P^*:\cO(Y_1) \ra A_0$, is generated
by $T$; we then   say that $T$ are {\it uniform} coordinates on $Y_1$. It then follows that
 the sheaves $\Omega^i_{Y^{\widehat{P}}}$ are free and
$$H^0(Y^{\widehat{P}},\Omega^i_{Y^{\widehat{P}}})=\bigoplus_{j_1<...<j_i} A_0[[T]]dT_{j_1}\wedge ... \wedge dT_{j_i}.$$

 {\it 3) Case} $Z=Z' \cap Z''$ where the ideal of $Z'$ is generated by $p_k$ for some $k$ and $Z''$ is the image in $Y$ of a  point $P\in Y(A_0)$; hence $Z$ is defined by the ideal $I=(p_k,P)$. In this case
$$H^0(Y^{\widehat{(p_k,P)}},\Omega^i_{Y^{\widehat{(p_k,P)}}})=\bigoplus_{j_1<...<j_i} \bZ_{p_k}[[T]]dT_{j_1}\wedge ...\wedge dT_{j_i}.$$

\medskip

 Corresponding to examples 1) and 2) above one can introduce more terminology as follows.
 Assume $Y$ is a smooth affine scheme over $A_0$ with irreducible geometric fibers, let $P\in Y(A_0)$ be a point, and let $\omega$ be a $1$-form on $Y$.

1) If  $\omega$ is  $\cP$-adically exact then there is a unique $\cP$-adic primitive
\begin{equation}
\label{fammm}
f=(f_k)_k=(f_1,...,f_d)\in \cO(Y^{\widehat{p_1...p_d}})=\cO(Y^{\widehat{p_1}}) \times ... \times \cO(Y^{\widehat{p_d}})\end{equation}
 of $\omega$ with $f_k(P)=0$ for all $k=1,...,d$.
(Here $f_k(P)$ is the image of $f_k$ via the map $\cO(Y^{\widehat{p_k}}) \ra \bZ_{p_k}$ defined by $P$.)
We shall refer to $f=(f_k)_k$ as the $\cP$-{\it adic primitive of $\omega$ normalized along $P$}.

2) If $\omega$ is $P$-adically exact then there is a unique $P$-adic primitive $f_0 \in \cO(Y^{\widehat{P}})$
of $\omega$ such that $f_0(P)=0$. (Here
 $f_0(P)$ is the image of $f_0$ via the map $\cO(Y^{\widehat{P}})\ra A_0$ defined by $P$.)
We shall refer to  $f_0$ as the $P$-{\it adic primitive of $\omega$ normalized along $P$}.

\medskip

\begin{remark}
The concepts above are directly related to the concept of {\it analytic continuation between primes} considered in \cite{laplace}.
Indeed
assume $Y$ is a smooth affine scheme over $A_0$ with irreducible geometric fibers, let $P\in Y(A_0)$ be a uniform point, and let $\omega$ be a $1$-form on $Y$.
Assume that $\omega$ is both $\cP$-adically exact and $P$-adically exact.
(This will be often the case in the main examples to be encountered later in the paper.)
Let $f=(f_k)$ and  $f_0$ be the $\cP$-adic and the $P$-adic primitives of $\omega$ normalized along $P$ respectively.
Since $df_k=df_0$ in $\bigoplus_j \bZ_{p_k}[[T]]dT_j$ and $f_k(P)=f_0(P)$ in $\bZ_{p_k}$ it follows that
 $f_k=f_0$ in $\bZ_{p_k}[[T]]$ which shows that $f$ is {\it analytically continued along $P$} and is {\it represented by $f_0$}
 in the sense of \cite{laplace}, Definition 2.23.
\end{remark}

\begin{remark}
The concepts above can be used to define a certain type of integration (and periods) of $1$-forms
in the arithmetic setting as follows. Indeed
assume again that $Y$ is a smooth affine scheme over $A_0$ with irreducible geometric fibers and let $\cS \subset Y(A_0)$ be a non-empty set of  points. (In our applications $\cS$ will be the whole of $Y(A_0)$.)
By an {\it elementary chain} we will understand a pair $(P_1,P_2)$ of points $P_1\in Y(\bZ_{p_{k_1}})$,
$P_2 \in Y(\bZ_{p_{k_2}})$ such that either 1) $k_1=k_2=:k$ or 2) $k_1\neq k_2$ and there exists $P_{12}\in \cS$ inducing both $P_1$ and $P_2$.
In case 1) we say $(P_1,P_2)$ is {\it vertical} and we let $\overline{P_1P_2}$ denote the ideal $(p_k)$ in $\cO(Y)$. In case 2) we say $(P_1,P_2)$ is {\it horizontal} and denote by
$\overline{P_1P_2}$ the ideal $P_{12}\subset\cO(Y)$.
By a {\it chain} we understand a tuple $\Gamma:=(P_1,...,P_N)$
where $(P_1,P_2),...,(P_{N-1},P_N)$ are elementary chains.
A chain as above is called a {\it cycle} if $P_1=P_N$.
We define the {\it group of abstract periods}
	$$
	\Pi:=\frac{\bigoplus_{k=1}^d \bZ_{p_k}}{\{(a_1,...,a_d)\in A_0^d;\sum_{k=1}^d a_k=0\}}.
	$$
For each $k$ the natural morphism  $\bZ_{p_k} \ra \Pi$ is injective,
and we shall view it as an inclusion.
Now let $\omega$ be a $1$-form on $Y$ which
is both $\cP$-adically exact and $\cS$-adically exact. Then one can
define the {\it integral of $\omega$ along a chain $\Gamma=(P_1,...,P_N)$} by the formula
\begin{equation}
\label{integral}
\int_{\Gamma} \omega:=\sum_{j=1}^{N-1} (f_j(P_{j+1})-f_j(P_j))\in \Pi,\end{equation}
where $f_j$ is any $\overline{P_jP_{j+1}}$-adic primitive of $\omega$.
(Here, if $(P_j,P_{j+1})$ is horizontal, then $f_j(P_j)$ and $f_j(P_{j+1})$ are defined as the
corresponding images of $f_j$ via
the homomorphisms $\cO(Y^{\widehat{P_{12}}})\ra A_0 \subset \bZ_{p_{k_j}}$  and
$\cO(Y^{\widehat{P_{12}}})\ra A_0 \subset \bZ_{p_{k_{j+1}}}$
defined by $P_{12}$ so they are equal; hence, in the sum (\ref{integral})
the terms corresponding to horizontal elementary chains are equal to $0$.
Also, note that summation above gives a well-defined element of $\bigoplus_k \bZ_{p_k}$, not
only of $\Pi$.) If $\Gamma$ is a cycle we may
refer to $\int_{\Gamma} \omega$ as a {\it period}.
Note that if $\omega$ is exact on $Y$ then its periods vanish.
(Indeed, for $\omega$ exact on $Y$  the summation (\ref{integral}) viewed as an element of $\bigoplus_k \bZ_{p_k}$, although generally not zero, becomes zero in $\Pi$.)

In the main examples to be later encountered in the
paper  the periods will be typically non-zero. In contrast with this phenomenon, in the projective (as opposed to the affine) case, and under certain Arakelov-style conditions at infinity, forms that are $\cP$-adically and $\cS$-adically exact should be expected to be exact; cf. \cite{bost}, where arithmetic analogues of the Hironaka-Matsumura theorems \cite{HM} from formal geometry are proved.
\end{remark}

\section{Arithmetic jet spaces: review  and complements}
In this section we review some of the concepts introduced in \cite{char, book, laplace, borger} and
provide some complements to that theory.

Let $C_p(X,Y) \in \bZ[X,Y]$ be the polynomial with
integer coefficients
	\[
	C_p(X,Y):=\frac{X^p+Y^p-(X+Y)^p}{p}.
	\]
A $p-${\it derivation} from a ring $A$ into an $A-$algebra
$\varphi:A \ra B$ is a set map $\d:A \ra B$ such that
	\begin{align}
		\d(1) &= 0, \label{delta-ring-axiom-1}\\
		\d(x+y) & = \d x + \d y +C_p(x,y), \label{delta-ring-axiom-2}\\
		\d(xy) & = \varphi(x^p) \cdot \d y + \varphi(y^p) \cdot\d x +p \cdot \d x \cdot \d y, \label{delta-ring-axiom-3}
	\end{align}
for all $x,y \in A$. Given a
$p-$derivation we always denote by $\phi:A \ra B$ the map
	\begin{equation} \label{eq:frob-lift}
		\phi(x)=\varphi(x)^p+p \d x;
	\end{equation}
then $\phi$ is a ring homomorphism.
Conversely, if $p$ is not a zero divisor in $B$
and if $\phi:A\to B$ is a ring homomorphism satisfying
the Frobenius lift property $\phi(x)\equiv \varphi(x)^p\bmod pB$
for all $x\in A$,
then there is a unique $p$-derivation $\delta:A\to B$ satisfying (\ref{eq:frob-lift}).

For any two distinct rational primes $p_1,p_2$ consider the
polynomial $C_{p_1,p_2}$ in the ring $\bZ[X_0,X_1,X_2]$ defined by
\begin{equation}
\label{commutator}
C_{p_1,p_2}(X_0,X_1,X_2) := \frac{C_{p_2}(X_0^{p_1},p_1X_1)}{p_1}
-\frac{C_{p_1}(X_0^{p_2},p_2 X_2)}{p_2}
 -\frac{\delta_{p_1} p_2}{p_2} X_2^{p_1}+ \frac{\delta_{p_2}p_1}{p_1}
X_1^{p_2}\, .
\end{equation}
Let $\cP=\{p_1,\ldots ,p_d\}$ be a finite set of  primes in $\bZ$.
A $\d_{\cP}$-{\it ring} is a ring $A$  equipped with
$p_k$-derivations $\delta_{p_k}:A \ra A$,  $k=1,\ldots ,d$, such that
\begin{equation}
\label{identit}
\delta_{p_k}\delta_{p_l}a-\delta_{p_l}\delta_{p_k}a=C_{p_k,p_l}(a,
\delta_{p_k}a,\delta_{p_l}a)
\end{equation}
for all $a \in A$, $k,l=1,\ldots ,d$.
A {\it homomorphism
of $\delta_{\cP}$-rings} $A$ and $B$ is a homomorphism
of rings $\varphi: A \rightarrow B$ that commutes
with the $p_k$-derivations in $A$ and $B$, respectively.
If $\phi_{p_k}(x)=x^{p_k}+p_k\d_{p_k}x$ is the homomorphism  associated to $\delta_{p_k}$,
condition
(\ref{identit}) implies that
\begin{equation}
\label{commu}
\phi_{p_k}\phi_{p_l}(a)=\phi_{p_l}\phi_{p_k}(a)
\end{equation}
for all $a \in A$. Conversely, if the commutation relations (\ref{commu})
hold, and the numbers $p_k$ are not zero divisors in $A$, then conditions
(\ref{identit}) hold, and we have that
$
\phi_{p_k}\delta_{p_l} a  =  \delta_{p_l}\phi_{p_k}a
$
for all $a \in A$. If $A$ is a $\d_{\cP}$-ring then for all $k$, the $p_k$-adic completions
$A^{\widehat{p_k}}$ are $\d_{\cP}$-rings in a natural way.

For a relation between these concepts and the theory of $\lambda$-rings we refer to
\cite{borger} and the references therein.

We let $\bZg=\{0,1,2,3, \ldots \}$, and let $\bZg^d$ be given the product order. We let $e_k$ be the
element of $\bZg^d$ all of whose components are zero except the $k$-th, which is $1$. We set $e=\sum
e_k=(1,\dots,1)$. For $i=(i_1,\dots,i_d) \in \bZg^d$ we set $\cP^i=p_1^{i_1}\cdots p_d^{i_d}$ and
$\d_{\cP}^i=\d_{p_1}^{i_1}\circ\cdots\circ\d_{p_d}^{i_d}$ and
$\phi_{\cP}^i=\phi_{\cP^i}=\phi_{p_1}^{i_1}\circ\cdots\circ\phi_{p_d}^{i_d}$.
A {\it $\delta_{\cP}$-prolongation system} $A^*=(A^r)$ is an inductive
system of rings $A^r$ indexed by $r \in \bZg^d$, provided with
transition maps $\varphi_{rr'}:A^r \ra A^{r'}$ for any pair of indices
$r$, $r^{'}$ such that $r \leq r'$, and equipped with $p_k$-derivations
$$
\delta_{p_k}:A^r\ra A^{r+e_k}\, ,
$$
$k=1,\ldots ,d$, such that (\ref{identit}) holds for all $k$, $l$, and
such that
$$
\varphi_{r+e_k,r'+e_k}\circ \delta_{p_k}=\delta_{p_k} \circ
\varphi_{rr'} :A^r \ra A^{r'+e_k}
$$
for  all $r\leq r'$ and all $k$.
A morphism of prolongation systems $A^* \ra B^*$ is a system of ring
homomorphisms
$u^r:A^r \ra B^r$ that commute with all the maps
$\varphi_{rr'}$ and the $\delta_{p_k}$ of $A^*$ and $B^*$.

Any $\delta_{\cP}$-ring $A$ induces a $\delta_{\cP}$-prolongation system $A^*$ where $A^r=A$ for all $r$
and $\varphi$ is the identity map.
Conversely, if $(A^r)$ is a $\d_{\cP}$-prolongation system then the ring
$$A^{\infty}:=\lim_{\ra} A^r$$ has a natural structure of $\d_{\cP}$-ring.

Let us say that a $\d_{\cP}$-prolongation system $(A^r)$ is {\it faithful} if the primes in $\cP$ are
non-zero divisors in all the rings $A^r$ and if all the homomorphisms $\varphi_{rr'}:A^r \ra A^{r'}$ are
injective. When this is the case we will usually view $\varphi_{rr'}$ as inclusions $A^r \subset A^{r'}$
and the primes in $\cP$ are non-zero divisors in $A^{\infty}$.

If $A^*=(A^r)$ is a $\d_{\cP}$-prolongation system then for each $p_k \in \cP$ the system of $p_k$-adic
completions $((A^r)^{\widehat{p_k}})_r$ is easily seen to have a unique structure of
$\d_{\cP}$-prolongation system with the property that
the natural maps $A^r \ra (A^r)^{\widehat{p_k}}$ define a morphism of prolongation systems.
Further, the operators $\d_{p_k}$ are continuous.

Let us say that a $\d_{\cP}$-ring  $A$ is $\d_{\cP}$-{\it generated} by a subring $A^0\subset A$ if
$A$ is generated as an $A^0$-algebra by the set
$\{\d^s_{\cP} a;\ a\in A^0,\ s \geq 0\}.$

For any affine scheme of finite type $X$ over $A_0=S^{-1}\bZ$ one can define a system of schemes of
finite type, $\cJ^r_{\cP}(X)$ over $A_0$, called the {\it $\d_{\cP}$-jet spaces} of $X$
(or {\it $\cP$-typical jet spaces} in the language of \cite{borger}); if $X=\Spec
A^0$, with $A^0=A_0[T]/(f)$, $T$ a tuple of variables, and $f$ a tuple of polynomials, then $\cJ^r_{\cP}(X)=\Spec A^r$, where
	\begin{equation} \label{eq:presentation-of-jet-space}
		A^r=A_0[\d_{\cP}^i T;i\leq r]/(\d_{\cP}^if;i\leq r).
	\end{equation}
Here, each $\d_{\cP}^i T$ is to be interpreted as a tuple of free variables, and each relation $\d_{\cP}^i f$
is interpreted as a polynomial expression in the variables $\d_{\cP}^i T$ by expanding it using
the rules (\ref{delta-ring-axiom-1})--(\ref{delta-ring-axiom-3}). The family $(A^r)_r$ forms
a $\d_{\cP}$-prolongation system by $\d_{p_k}(\d_{\cP}^i T):=\d_{\cP}^{i+e_k} T$.
In the case where $\cP$ consists of a single prime $p$, formal $p$-adic versions of these
spaces were introduced in \cite{char}. In the multiple-prime case considered here, they were introduced
(independently) in \cite{laplace} and \cite{borger}.
In \cite{borger} the spaces $\cJ_{\cP}^r(X)$ were denoted by $W_{r*}(X)$; the notation here follows
\cite{laplace}.

The prolongation system $(A^r)$, where $A^r$ is still $\cO\big(\cJ_{\cP}^r(X)\big)$, has
the following universal property: for any $\d_{\cP}$-prolongation system $(B^r)$, each
ring map $A^0\to B^0$ extends uniquely to a map $(A^r)\to(B^r)$ of $\d_{\cP}$-prolongation systems.
Let us now apply this in the case where $B^r$ is the constant inductive system $A_0$, where each $B^r$ is
$A_0$ and each $\varphi_{rr'}$ the identity map. Since $A_0$ has a unique
$\d_{\cP}$-ring structure, this system has a unique $\d_{\cP}$-prolongation structure.
By the universal property, any ring map $A^0\to A_0$ extends to a unique map $A^r\to B^r=A_0$
of $\d_{\cP}$-prolongation systems, and hence defines a map
$X(A_0)\to \cJ_{\cP}^r(X)(A_0)$ for each $r$.  We call the image of a point $P\in X(A_0)$ under this
map its \emph{canonical lift} (to $\cJ_{\cP}^r(X)$). We will often denote it by $P^r$ (because of
its relation with the map (\ref{eq:coghost-map}) below).

For $A_0$-algebras $C$, there is a functorial isomorphism
	\begin{equation} \label{eq:Witt-jet-adjunction}
		\Hom_{A_0}(\Spec C, \cJ_{\cP}^r(X)) = \Hom_{A_0}(\Spec W_r(C), X)
	\end{equation}
Here $W_r(C)$ denotes the $A_0$-algebra of $\cP$-typical Witt vectors of length $r$ with entries in $C$.
If $\cP$ consists of a single prime $p$, then $W_r(C)$ is the usual ring of $p$-typical Witt vectors of
length $r$ (length $r+1$ in the traditionally more common numbering) with entries in $C$. In
\cite{borger}, the equation (\ref{eq:Witt-jet-adjunction}) is taken to be the definition of
$\cJ_{\cP}^r(X)$, or rather of its functor of points. This lets us define $\cJ_{\cP}^r(X)$ for
rather general $X$, such as algebraic spaces.
We will not need this generality here, but we will make use of (\ref{eq:Witt-jet-adjunction}), as well
as some results in \cite{borger} proved using this functorial point of view, such as the
following:

\begin{proposition}\label{pro:jets-preserve-smoothness}
	Let $f:X'\ra X$ be a smooth (respectively \'etale)
	morphism. Then the induced morphism $g:\cJ_{\cP}^r(X')\ra \cJ_{\cP}^r(X)$ is smooth (respectively
	\'{e}tale). If $f$ is also surjective, then so is $g$.
\end{proposition}

The rings $A^r=\cO(\cJ^r_{\cP}(X))$ form a $\d_{\cP}$-prolongation system with $A^{\infty}$ being
$\d_{\cP}$-generated by $A^0$.
If $Y \subset X$ is a principal open set of $X$, $\cO(Y)=\cO(X)_f$ then $\cO({\mathcal
J}_{\cP}^r(Y)) \simeq \cO({\mathcal J}_{\cP}^r(X))_{f_r}$ where $f_r=\prod_{i\leq r}\phi_{\cP}^i(f)$. In
particular, the induced morphism ${\mathcal J}_{\cP}^r(Y) \ra {\mathcal J}_{\cP}^r(X)$ is an open
immersion.

In corollary \ref{rece}, we will show that if $X/A_0$ smooth and $\cJ_{\cP}^r(X)=\Spec A^r$
then the prolongation system $(A^r)$ is faithful.

Let us note that the homomorphisms
$$\varphi_{sr} \circ \phi_{\cP}^s:A^0\ra A^r$$
for $s \leq r$ induce morphisms $\cJ^r_{\cP}(X)\ra X$, and hence
a natural morphism
	\begin{equation} \label{eq:coghost-map}
		\cJ^r_{\cP}(X)\longlabelmap{\kappa_{\leq r}} \prod_{s \leq r} X.
	\end{equation}
This map becomes an isomorphism after tensoring with $\bZ[1/p_1,\dots,1/p_k]$.

\begin{proposition}\label{pro:jet-group-maps}
	If $X$ is a group scheme, then the maps $\cJ^s_{\cP}(X) \to \cJ^r_{\cP}(X)$
	and $\cJ^s_{\cP}(X) \to X$ induced by $\varphi_{sr}$ and $\phi^s_{\cP}$ are group homomorphisms.
\end{proposition}

\begin{proof}
It is enough to show that each map induces a group homomorphism on $B$-valued points for each
$A_0$-algebra $B$. We can describe these maps simply using (\ref{eq:Witt-jet-adjunction}).
For $\varphi_{sr}$, it can be identified with the map
	$$
	X(W_s(B))\to X(W_r(B))
	$$
induced by a map $W_s(B)\to W_r(B)$ (namely, the natural projection map).
For $\phi^s_{\cP}$, it can be identified with the map
	$$
	X(W_s(B)) \to X(B)
	$$
induced by a map $W_s(B)\to B$ (namely, the $s$-th ghost component map).
In particular, each map is a group homomorphism.
\end{proof}

%

\begin{remark}
\label{nonaffinejetspaces}
Let $X$ be a smooth scheme over $A_0$ and consider an affine open cover $X=\bigcup X_i$.
Let $p$ be a prime and let $r \geq 0$ be an integer.
Then one can glue the formal schemes $\cJ_{\{p\}}^r(X_i)^{\widehat{p}}$ along the formal schemes $\cJ_{\{p\}}^r(X_i \cap X_j)^{\widehat{p}}$
to construct a formal scheme that we denote by $J^r_p(X)$.  This latter formal scheme does not depend on the covering we chose for $X$ and will be referred to as the $p$-{\it jet space} of $X$ of order $r$.  After base change to  $\hat{\bZ}_p^{ur}$ (the completion of the maximum unramified extension of $\bZ_p$), $J^r_p(X)$ becomes  equal to the  $p$-jet space of $X$  considered  in \cite{char, book}.
Note that taking formal completions here is needed for this approach to work. Otherwise,
one would use the method of \cite{borger}.
\end{remark}

\begin{proposition}\label{pro:smooth-transition-maps}
	Let $X$ be a smooth affine scheme
	over $A_0$, and let $r,s\in\bZ_{\geq 0}^d$ be elements
	with $r\leq s$.  Then the map $\cJ^s_{\cP}(X)\to\cJ^r_{\cP}(X)$ induced by $\varphi_{sr}$ is
	smooth and surjective.
\end{proposition}

\begin{proof}
	First, let us reduce to the case where the number $d$ of primes in $\cP$ is $1$.
	The case $d=0$ is trivial.  Now assume the theorem when $d=1$,
	and suppose $d\geq 2$.  Let $\cP'$
	denote the set $\{p_1,\dots,p_{d-1}\}$, let
	$s'$ denote $(s_1,\dots,s_{d-1})$, and let $r'$
	denote $(r_1,\dots,r_{d-1})$.  Then we have the following factorization of
	the map in question:
		$$
		\xymatrix{
		\cJ^s_{\cP}(X) \ar[rr] \ar^{\sim}[d]
			&
			& \cJ^r_{\cP}(X) \ar^{\sim}[d]\\
		\cJ^{s_d}_{\{p_d\}}\cJ^{s'}_{\cP'}(X) \ar^a[r]
			& \cJ^{s_d}_{\{p_d\}}\cJ^{r'}_{\cP'}(X) \ar^b[r]
			& \cJ^{r_d}_{\{p_d\}}\cJ^{r'}_{\cP'}(X).
		}
		$$
	By proposition \ref{pro:jets-preserve-smoothness},
	the space $\cJ^{r'}_{\cP'}(X)$ is smooth over $\cJ^{r'}_{\cP'}(\Spec A_0) = \Spec A_0$.
	Therefore
	by the case $d=1$, the map $b$ is smooth and surjective.  And by induction,
	the map $\cJ^{s'}_{\cP'}(X)\to \cJ^{r'}_{\cP'}(X)$ is smooth and surjective.
	By~\cite{borger} (7.2(c)),
	the map $a$ is then smooth and surjective, and therefore $b\circ a$ is.

	Thus it suffices to assume $d=1$.  Let us write $p=p_1$, $\cJ=\cJ_{\{p\}}$ and so on.
	We can further assume $s=r+1$ with $r\geq 0$, again by induction.
	Let $\varphi$ denote the map $\cJ^{r+1}(X)\to\cJ^r(X)$.
	
	To show $\varphi$ is smooth, we need to show it is formally smooth and locally of finite
	presentation. Finite presentation is clear: by (\ref{eq:presentation-of-jet-space}),
	both $\cJ^{r+1}(X)$	and $\cJ^r(X)$ are of finite type over $A_0$, which is Noetherian.
	
	Let us now show that $\varphi$ is formally smooth.
	Let $g:B\to A$ be a surjection of rings with square-zero kernel.  By definition,
	we need to show that the evident map
		$$
		\cJ^{r+1}(X)(B) \longrightarrow \cJ^{r+1}(X)(A) \times_{\cJ^r(X)(A)} \cJ^r(X)(B)
		$$
	is a surjection.  By (\ref{eq:Witt-jet-adjunction}),
	we can identify this map with the evident map
		$$
		X(W_{r+1}(B)) \longrightarrow X(W_{r+1}(A)) \times_{X(W_r(A))} X(W_r(B)).
		$$
	Since $X$ is affine, this map can be identified with the evident map
		$$
		X(W_{r+1}(B)) \longrightarrow X(W_{r+1}(A)\times_{W_r(A)} W_r(B)).
		$$
	Thus, since $X$ is formally smooth over $\bZ$, we only need to show that the evident
	ring map
		\begin{equation} \label{eq:smooth-jet-1}
			W_{r+1}(B) \longrightarrow W_{r+1}(A)\times_{W_r(A)} W_r(B)
		\end{equation}
	is surjective with nilpotent kernel.
	Using the usual Witt components, we can identify this with the map
		$$
		B^{r+2} \longrightarrow A^{r+2}\times_{A^{r+1}} B^{r+1},
		$$
	defined by
		$$
		(b_0,\dots,b_{r+1}) \mapsto \Big(\big(g(b_0),\dots,g(b_{r+1})\big),(b_0,\dots,b_r)\Big).
		$$
	This map can further be identified with the map
		\begin{equation} \label{eq:smooth-jet-2}
		B^{r+2} \longrightarrow B^{r+1} \times A
		\end{equation}
	defined by $(b_0,\dots,b_{r+1})\mapsto \big(b_0,\dots,b_r,g(b_{r+1})\big)$.
	Since $g$ is surjective, so is (\ref{eq:smooth-jet-2}), and hence so is (\ref{eq:smooth-jet-1}).
	On the other hand, the kernel of (\ref{eq:smooth-jet-2}) is the set of
	elements $(0,\dots,0,b)$, where $b\in\ker(g)$.
	This kernel then has square zero, since we have
		\begin{equation}\label{eq:witt-mult}
		(0,\dots,0,b)(0,\dots,b') = (0,\dots,0,p^{r+1}bb') = 0
		\end{equation}
	and $\ker(g)^2=0$. Therefore (\ref{eq:smooth-jet-1}) has square-zero kernel.
	It follows that $\varphi$ is formally smooth.
	
	Let us finally show that $\cJ^{r+1}(X)\to\cJ^r(X)$ is surjective.
	It suffices to show this after base change to $\bZ[1/p] \times\bF_p$.
	Over $\bZ[1/p]$, the map can be identified with the projection
	$X^{r+2}\to X^{r+1}$, which is surjective.
	
	For the base change to $\bF_p$,
	we will show the stronger property that for any $\bF_p$-algebra $A$, the map
		$$
		\cJ^{r+1}(X)(A) \longrightarrow \cJ^r(X)(A)
		$$
	is surjective.  By (\ref{eq:Witt-jet-adjunction}), this can be identified with the map
		$$
		X(W_{r+1}(A)) \longrightarrow X(W_r(A)).
		$$
	But this map is surjective because $X$ is formally smooth and because
	the map $W_{r+1}(A)\to W_r(A)$ is
	a surjection with nilpotent kernel, again by (\ref{eq:witt-mult}).
\end{proof}

\begin{remark}
Proposition \ref{pro:smooth-transition-maps} holds, more generally, for $X$ an arbitrary smooth
algebraic space over $A_0$ (See
\cite{borger} for the definition of arithmetic jet spaces of algebraic spaces.)
The proof is as follows: by general \'etale localization properties of $\cJ_{\cP}^r$, it is enough
to replace $X$ with an affine \'etale cover, and in this case the result was shown above.
We will not need the non-affine version in this paper.
\end{remark}

\begin{corollary}
\label{rece}
If $\cJ_{\cP}^r(X)=\Spec A^r$ then  the $\d_{\cP}$-prolongation systems
$(A^r)$ and $((A^r)^{\widehat{p_k}})$ are faithful.
\end{corollary}

{\it Proof}.
By proposition \ref{pro:smooth-transition-maps} the transition maps $A^s \ra A^r$ are faithfully flat and
hence injective. Also since $A^0$ is flat over $\bZ$ the primes in $\cP$ are not zero divisors in any of
the rings $A^r$. Hence $(A^r)$ is faithful. Now the primes in $\cP$ continue to be non-zero divisors in
$(A^r)^{\widehat{p_k}}$. Moreover the maps $A^s/p_kA^s \ra A^r/p_kA^r$ are faithfully flat hence
injective. It follows that the maps $(A^s)^{\widehat{p_k}}\ra(A^r)^{\widehat{p_k}}$ are injective. So
$((A^r)^{\widehat{p_k}})$ is faithful. \qed

\bigskip

We also have the following useful result:

\begin{proposition}
\label{codim2}
Let $X$ be a smooth affine scheme over $A_0$ and let $X=\bigcup_i X_i$ be a covering with principal affine open sets. Then, for each $r$,  the closed set
\begin{equation}
\label{diffe}
\cJ^r_{\cP}(X) \setminus \bigcup_i \cJ_{\cP}^r(X_i) \end{equation}
has codimension $\geq 2$ in $\cJ_{\cP}^r(X)$.
\end{proposition}

{\it Proof}. Let $U=\bigcup_i \cJ_{\cP}^r(X_i)$ and $Y=\cJ^r_{\cP}(X)$. It is enough to show
that for any prime $p$, the complement $(Y\otimes \bF_p) \setminus (U \otimes \bF_p)$ has codimension
$\geq 2$ in $Y \otimes \bF_p$ and that $(Y\otimes \bQ) \setminus (U \otimes \bQ)$
has codimension $\geq 2$ in $Y \otimes \bQ$.
This follows exactly as in the proof of Proposition 2.22 in \cite{laplace}.
\qed

\begin{remark}
We will later need the following consequence of \cite{borger}, Proposition 8.2 (cf. \cite{laplace},
Remark 2.27). Let $X$ be a smooth affine scheme over $A_0$ with connected geometric fibers and let
$P\in X(A_0)$ be a uniform point with uniform coordinates $T$ in $X$. Let $Y=\cJ_{\cP}^r(X)$ and let
$P^r \in Y(A_0)$ be the canonical lift of $P$, as defined above.
Then $P^r$ is a uniform point of $Y$ with uniform coordinates $(\d_{\cP}^iT;i\leq r)$ in $Y$.
\end{remark}

\section{The tangent bundle of arithmetic jet spaces}

In this section we extend the theory in \cite{book}, pp. 88-100, from the case where $\cP$ consists
of one prime to the case where it consists of several primes.

Let $A$ be a $\d_{\cP}$-ring in which the primes in $\cP$ are non-zero divisors and let $A^0\subset A$ be
a subring. Let $\partial:A^0\ra A^0$ be a derivation and let $r \in \bZ_{\geq 0}^d$. A derivation
$\partial_r:A \ra A$ will be called an {\it $r$-conjugate} of $\partial$ on $A$ if for any $s\in
\bZ_{\geq 0}^d$ we have
	\begin{equation} \label{starstar}
		\partial_{r} \circ \phi_{\cP}^{s}
			=\delta_{rs} \cdot \cP^r \cdot \phi^{s}_{\cP} \circ \partial:A^0\ra A,
	\end{equation}
where $\d_{rs}$ is the Kronecker symbol.

In other words, let $X=\Spec A^0$, let $\kappa:\cJ_{\cP}(X) \longlabelmap{} \prod_r X$ denote the
limit of the maps $\kappa_{\leq r}$ of (\ref{eq:coghost-map}),
let $\partial'$ denote the vector field $(\dots,\cP^r\partial,\dots)_r$ on $\prod_r X$, and
let $\partial''$ denote the vector field on $\prod_r X$ with values in $A$ induced by $\partial'$.
Then an $r$-conjugate of $\partial$ is an extension of $\partial''$ to
a vector field $\partial_r$ on $X=\Spec A$.

A system of derivations $(\partial_{r})_{r}$, $\partial_{r}:A \ra A$,
indexed by multi-indices $r\in \bZ_{\geq 0}^d$, will be called a {\it complete system of
conjugates} of $\partial$ on $A$ if for any $r$ the derivation $\partial_{r}$ is an $r$-conjugate of
$\partial$ on $A$.

Clearly we have the following uniqueness result:

\begin{proposition}
\label{uniqueness}
Assume $A$ is a $\d_{\cP}$-ring in which the primes in $\cP$ are non-zero divisors
and assume $A$ is $\d_{\cP}$-generated by a subring $A^0 \subset A$.
Let $\partial:A^0\ra A^0$ be a derivation. Then, for each $r$, there is at most one $r$-conjugate
$\partial_r$ of $\partial$ in $A$. In addition we have $\partial_r A^n=0$ for $r \not\leq n$.
\end{proposition}

On the other hand we have the following existence result.

\begin{proposition}
\label{existence}
Let $X=\Spec A^0$ be an affine smooth scheme over $A_0$ and let $\cJ^r_{\cP}(X)=\Spec A^r$ be its
$\d_{\cP}$-jet spaces.
Let $\partial:A^0\ra A^0$ be a derivation. Then there exists a complete system of conjugates
$(\partial_r)_r$ of $\partial$ in $A^{\infty}$
(which is unique by Proposition \ref{uniqueness}). Moreover
if, by Corollary \ref{rece}, we view each $A^n$ as a subset of $A^{\infty}$,
then we have $\partial_r A^n \subset A^n$ for all $r,n$.
\end{proposition}

We denote by $\partial_{r|n}:A^n \ra A^n$ the restriction of $\partial_r$ to $A^n$.

\medskip

{\it Proof}.
Let $A^0=A_0[T]/(f)$ where $T$ is a tuple of indeterminates $T_a$ and $f$ is a tuple of polynomials $f_b$.
Lift $\partial$ to a derivation $\partial:A_0[T]\ra A_0[T]$ and, for any $n$ define the derivation $\partial_r$
on $\bQ[\d_{\cP}^{s}T;s\geq 0]=\bQ[\phi_{\cP}^{s}T;s\geq 0]$ as the unique derivation satisfying
\begin{equation}
\label{heart}
\partial_r(\phi_{\cP}^{s}T_a)=\d_{rs} \cdot \cP^r \cdot \phi_{\cP}^{s}(\partial T_a).
\end{equation}
Clearly (\ref{starstar}) holds and $\partial_r$ sends each $\bQ[\d_{\cP}^sT;s \leq n]$ into itself. Set $\partial_{r}=0$ for $r\in \bZ^d \setminus \bZ_{\geq 0}^d$.

\medskip

{\it Claim 1}. $\partial_r \circ \phi_{p_k}=p_k\cdot \phi_{p_k}\circ \partial_{r-e_k}$ on $\bQ[\d_{\cP}^{s}T;s\leq n]$.

\medskip

Indeed the difference between the left hand side and the right hand side of the above equality is a derivation that vanishes on the set of generators
$\{\phi_{\cP}^{s}T;s\leq n\}$ of $\bQ[\d_{\cP}^{s}T;s\leq n]$; for, using (\ref{heart}), we have
$$\begin{array}{rcl}
(\partial_r \phi_{p_k})(\phi_{\cP}^{s}T_a) & = & \partial_r \phi_{\cP}^{s+e_k} T_a\\
\  & \  & \  \\
\  & = & \d_{r,s+e_k} \cdot \cP^r \cdot \phi_{\cP}^{s+e_k} \partial T_a\\
\  & \  & \  \\
\  & = & \d_{r-e_k,s} \cdot \cP^r \cdot \phi_{p_k}(\phi_{\cP}^{s}(\partial T_a))\\
\  & \  & \  \\
\  & = & p_k \cdot \phi_{p_k}(\cP^{r-e_k} \cdot \d_{r-e_k,s} \cdot \phi_{\cP}^{s}(\partial T_a))\\
\  & \  & \  \\
\  & = & (p_k \cdot \phi_{p_k} \circ \partial_{r-e_k})(\phi_{\cP}^{s}T_a).
\end{array}$$

\medskip

{\it Claim 2}. $\partial_r$ maps $A_0[\d_{\cP}^{s}T;s\leq n]$ into itself.

\medskip

Indeed it is enough to show that $\partial_r(\d_{\cP}^{s}T_a) \in A_0[\d_{\cP}^{s}T;s\leq n]$ for all $s\leq n$ and all $a$. We proceed by induction
on the sum of the components of $r+s$. The statement is clearly true for $s=0$ so for $r+s=0$. Now assume $r+s$ arbitrary. We may assume
$s\neq 0$ so we may assume there exists $k$ such that $s_k \geq 1$. Then
$$\begin{array}{rcl}
\partial_r(\d_{\cP}^{s} T_a) & = & \partial_r(\d_{p_k}\d_{\cP}^{s-e_k} T_a)\\
\  & \  & \  \\
\  & = & \partial_r\left( \frac{\phi_{p_k}(\d_{\cP}^{s-e_k}T_a)-(\d_{\cP}^{s-e_k}T_a)^{p_k}}{p_k}\right)\\
\  & \  & \  \\
\  & = & \phi_{p_k} \partial_{r-e_k}(\d_{\cP}^{s-e_k} T_a)-(\d_{\cP}^{s-e_k} T_a)^{p_k-1} \partial_r(\d_{\cP}^{s-e_k}T_a),
\end{array}$$
by Claim 1. The latter belongs to $A_0[\d_{\cP}^{s}T;s\leq n]$ by the induction hypothesis.

\medskip

{\it Claim 3}. $\partial_r$ maps the ideal $(\d_{\cP}^{rs}f;s\leq n)$ of  $A_0[\d_{\cP}^{s}T;s\leq n]$ into itself.

\medskip

Indeed one can repeat the argument in the proof of Claim 2 with $T_a$ replaced by $f_b$.

\medskip

Now Claims 2 and 3 imply that $\partial_r$ induce derivations $A^{\infty}\ra A^{\infty}$ satisfying (\ref{starstar})
which ends our proof.
\qed

\begin{remark}
The proof above shows that we have the commutation relations:
\begin{equation}
\label{getgon}
\partial_r \circ \phi_{\cP}^s=\cP^s \cdot \phi_{\cP}^s \circ \partial_{r-s}:A^{\infty} \ra A^{\infty}\end{equation}
for all $r,s$. (Here, as usual, $\partial_{r-s}=0$ if $s \not\leq r$.)
\end{remark}

\begin{remark}
Under the assumption of Proposition \ref{uniqueness} fix $k$ and consider the derivation (still denoted by ) $\partial$ on $(A^0)^{\widehat{p_k}}$
induced by $\partial$. Then the system of conjugates $(\partial_r)_r$ on $A^{\infty}$ induces a complete system of conjugates of $\partial:(A^0)^{\widehat{p_k}}\ra (A^0)^{\widehat{p_k}}$
on
$$A^{\infty}_k:=\lim_{\ra}(A^r)^{\widehat{p_k}}$$ which we continue to denote by $(\partial_r)_r$. Recall from Corollary \ref{rece} that
$(A^n)^{\widehat{p_k}}\subset A^{\infty}_k$ for all $n$.  Note that  we have
$$\partial_r((A^n)^{\widehat{p_k}}) \subset (A^n)^{\widehat{p_k}}.$$
\end{remark}

\begin{remark}
\label{ulle}
Under the assumptions of Proposition \ref{existence} we have the following formula:
$$(a \cdot \partial)_{r}=\phi_{\cP}^r(a) \cdot \partial_{r},$$
for any $a \in A$. This follows from the uniqueness in Proposition \ref{uniqueness}.
\end{remark}

\begin{proposition}
\label{collo}
Let $X=\Spec A^0$ be smooth over $A_0$ and assume $(\partial^a)_{1\leq a \leq m}$ is a basis for the tangent sheaf $T_{X}$ of $X$ over $A_0$. Then
the family of conjugates $(\partial_{r|n}^a)_{1 \leq a \leq m,\ 0\leq r \leq n}$
is a basis for the tangent sheaf $T_{\cJ^n_{\cP}(X)}$.
\end{proposition}

{\it Proof}.
Assume first $X$ has \'{e}tale coordinates, i.e. there is
  an \'{e}tale map $X \ra \bA^m=\Spec \ A_0[T]$, where $T$ is a family of indeterminates $T_b$. Note that, by (\ref{starstar}),
$$\begin{array}{rcl}
\cP^r\cdot \phi^r_{\cP}(\partial^a T_b) & = & \partial^a_{r}(\phi_{\cP}^r T_b)\\
\  & \  & \  \\
\  & = & \partial_{r}^a(\cP^r \cdot \d_{\cP}^r T_b+(\text{polynomial in $\d_{\cP}^s T$ with $s \leq r$, $s\neq r$}))\\
\  & \  & \  \\
\  & = & \partial_{r}^a(\cP^r \cdot \d_{\cP}^r T_b)\\
\  & \  & \  \\
\ & = & \cP^r \cdot \partial^a_{r}(\d_{\cP}^r T_b)\end{array}$$
\noindent Hence
\begin{equation}
\label{hii}\partial^a_{r}(\d_{\cP}^rT_b)= \phi_{\cP}^r(\partial^a T_b).\end{equation}
Similarly we get
\begin{equation}
\label{haa}
\partial^a_{r}(\d_{\cP}^s T_b)=0\end{equation}
for $s\leq r$, $s\neq r$.
Since $\cJ^n_{\cP}(X) \ra \cJ_{\cP}^n(\bA^m)$ is \'{e}tale we know that
$$\left( \frac{\partial}{\partial(\d_{\cP}^s T_b)} \right)_{1 \leq b \leq m,\ s \leq n}$$
is a basis of $T_{\cJ_{\cP}^n(X)}$.
So we may  write
$$\partial_{r}^a=\sum_{0\leq s \leq n}\sum_{b=1}^m
\alpha_{ra}^{sb} \cdot \frac{\partial}{\partial(\d_{\cP}^s T_b)},$$
with $\alpha_{ra}^{sb}\in A^n=\cO(\cJ^n_{\cP}(X))$. Now order the multiindices $(r,a)$ by the lexicographic order.
Then, by (\ref{hii}) and (\ref{haa}) the matrix $\alpha:=(\alpha_{ra}^{sb})$ consists of $m \times m$ blocks
such that all blocks under the diagonal are $0$ and the diagonal blocks are of the form $\phi_{\cP}^r \beta$,
where $\beta=(\partial^a T_b)$. Therefore the matrix $\alpha$ is invertible which ends our proof in case $X$ has \'{e}tale coordinates.

Assume now $X$ does not necessarily have \'{e}tale coordinates. Take a finite affine open cover $X=\cup X_i$ such that each
$X_i$ has \'{e}tale coordinates and let $X':=\coprod X_i$ be the disjoint union. Then $X'$ has \'{e}tale coordinates and is an \'{e}tale cover of $X$. By \cite{borger}, Propositions 7.2 and 7.4
the map $\pi^r:\cJ_{\cP}^r(X')\ra \cJ_{\cP}^r(X)$ is an \'{e}tale cover, in particular it is faithfully flat.
Consider the homomorphism $u:\cO_{\cJ_{\cP}^r(X)}^N\ra T_{\cJ_{\cP}^r(X)}$ defined by the collection of sections
$\partial_r^a$ of $T_{\cJ_{\cP}^r(X)}$ (where $N$ is the cardinality of this collection). The pull-back of $u$ to $\cJ_{\cP}^r(X')$ coincides with the homomorphism $u':\cO_{\cJ_{\cP}^r(X')}^N\ra T_{\cJ_{\cP}^r(X')}$ defined by the collection of corresponding sections (still denoted by)
$\partial_r^a$ of $T_{\cJ_{\cP}^r(X')}$; this follows from the uniqueness in Proposition \ref{uniqueness} plus the fact that
$\pi^r$ is \'{e}tale. By the first part of the proof
$u'$ is an isomorphism. Since  $\pi^r$ is faithfully flat it follows that $u$ itself is an isomorphism. This ends the proof.
\qed

\section{The cotangent bundle of arithmetic jet spaces}

\begin{lemma}
Let $X=\Spec  A^0$ be smooth over $A_0$ and let $\omega \in \Omega^i_{A^0}$ be a $i$-form on $X$.
Then $\phi^{r*}_{\cP}\omega \in \cP^{ir}\cdot \Omega^i_{A^r}$.
\end{lemma}

{\it Proof}. It is enough to show that $\phi_{p_k}^* \eta\in p_k \cdot \Omega_{A^{r+e_k}}$ for any
$\eta \in \Omega_{A^r}$. But if $\eta=fdg$ with $f,g\in A^r$ then
$$\phi_{p_k}^*\eta=(f^{p_k}+p_k\d_{p_k} f)d(g^{p_k}+p_k\d_{p_k}g)=(f^{p_k}+p_k\d_{p_k} f)(p_k g^{p_k-1}dg+p_kd(\d_{p_k}g)),$$
and we are done.
\qed

Since, for $X/A_0$ smooth, the primes in $\cP$ are non-zero divisors in all the rings
$A^r$ (and hence in $\Omega_{A^r}$), we may define, for any $i$-form
$\omega \in \Omega^i_{A^0}$,  the $i$-forms
	\begin{equation} \label{defome}
		\omega_r:=\frac{\phi_{\cP}^{r*}\omega}{\cP^{ir}}\in \Omega^i_{A^{r}}.
	\end{equation}
Furthermore, for $n \geq r$ we consider the forms
$\varphi_{rn}^*\omega_r\in \Omega^i_{A^n}$;
when $n$ is clear from the context we simply denote these forms, again, by  $\omega_r$.
Also, for $m:=\cP^r$ we write
\begin{equation}
\label{sqr}
\omega_{[m]}:=\omega_r.
\end{equation}
Note that if $\omega \in \Omega^i_{A^0}$ and $\eta \in \Omega^j_{A^0}$ then
\begin{equation}
\label{omegaeta}
(\omega \wedge \eta)_r=\omega_r \wedge \eta_r
\end{equation}
for all $r$.
Since the exterior derivative $d:\Omega_{A^r}^i\ra \Omega_{A^r}^{i+1}$ commutes with $\phi_{\cP}^{r*}$,
we have
	\begin{equation} \label{d-commute}
		(d\omega)_r = \cP^r d(\omega_r)
	\end{equation}
In particular, $\omega\in \Omega_{A^0}^i$ is closed then the forms $\omega_r \in \Omega_{A^r}^i$
are also closed for all $r$.

Finally, let us record the formula
	\begin{equation} \label{eq:divided-frob-formula}
		\frac{\phi_p^*}{p}(xdy) = \phi_p(x) \big(y^{p-1}dy + d(\delta_p y)\big),
	\end{equation}
for $p\in\cP$.

\begin{remark}
Of course, it is possible to take (\ref{eq:divided-frob-formula})
as the definition of $\phi_p^*/p$. We could then define
$\phi_{\cP}^{r*}/{\cP^r}$ as compositions of the single-prime operators $\phi_{p_k}^*/p_k$.
This approach has the benefit that it works in general, without any smoothness assumptions on $X$.
(See \cite{BW}, 12.5 and 12.8, where $\phi^*/p$ is denoted $\theta$.)
\end{remark}

\begin{remark}
Also these operators can be viewed as ``universal lifts" of the inverse Cartier operator on the corresponding
forms. To explain this assume for simplicity that $\cP=\{p\}$ and consider the case of $1$-forms. We recall the
definition of the inverse Cartier operator \cite{DI}. Assume $B$ is a smooth $\bF_p$-algebra. For any $B$-module
$M$ defined by $B \times M \ra M$, $(b,m)\mapsto bm$ we denote by $F_*M$ the additive group $M$ viewed as a
$B$-module under $B \times F_*M\ra F_*M$, $(b,m) \mapsto b \cdot m:=b^pm$. Then the deRham complex $$F_*
\Omega^*_B=(F_*B \stackrel{d}{\ra} F_*\Omega_B \stackrel{d}{\ra} F_*\Omega^2_B \ra \cdot \cdot \cdot )$$ is a
complex of $B$-modules, in particular $H^1(F_* \Omega^*_B)$ is a $B$-module. The map $B \ra H^1(F_*\Omega^*_B)$, $b
\mapsto [b^{p-1}db]$ (where $[\ ]$ means class in $H^1$) is a derivation so it induces a $B$-module homomorphism
$C^{-1}:\Omega_B \ra H^1(F_*\Omega^*_B)$, called the {\it inverse Cartier operator}; it satisfies $$C^{-1}(adb)=a
\cdot C^{-1}(db)=[a^p b^{p-1}db].$$

Now if $X=\Spec A^0$ and $A^r$ are as in the discussion before this remark and $\omega \in \Omega_{A^0}$ then
$\omega_1\in \Omega_{A^1}$. Let $\phi_0:(A^0)^{\widehat{p}} \ra (A^0)^{\widehat{p}}$ be any lift of Frobenius.
There is an induced section $s:X^{\widehat{p}} \ra J^1_p(X)$ of the projection $J^1_p(X)\ra X^{\widehat{p}}$ so
we may consider the pull-back $s^*\omega_1 \in \Omega_{X^{\widehat{p}}}$ which is of course nothing but
$s^*\omega_1=\frac{\phi_0^* \omega}{p}$. Denoting by upper bars reduction mod $p$ a trivial computation shows
that the class $[\overline{s^*\omega_1}]$ of $\overline{s^*\omega_1}\in \Omega_{A^0 \otimes \bF_p}$ in
$H^1(F_*\Omega^*_{A^0\otimes \bF_p})$ equals $C^{-1}(\overline{\omega})$, image of $\overline{\omega}\in
\Omega_{A^0\otimes \bF_p}$ under the inverse Cartier operator. (Cf. also \cite{DI} for the relation between the
Cartier operator and lifts of Frobenius.)
\end{remark}

\begin{proposition}
\label{morelli}
Let $X=\Spec A$ be smooth over $A_0$ of relative dimension $m$
 and assume $(\partial^a)_{1\leq a \leq m}$ is a basis for the tangent sheaf $T_{X}$ of $X$ over $A_0$.
Let $(\omega^a)_{1 \leq a \leq m}$ be the dual basis for $\Omega_{X}$ (i.e. $\langle  \omega^a, \partial^b \rangle=\d_{ab}$).
Then $(\omega^a_{r})_{1 \leq a \leq m,\ 0\leq r \leq n}$ is a basis for $\Omega_{A^n}$, dual to the basis $(\partial^a_{r|n})_{1 \leq a \leq m,0 \leq r \leq n}$.
\end{proposition}

{\it Proof}.
As in the proof of Proposition \ref{collo} we may assume
 there is an \'{e}tale map $X \ra \bA^m=\Spec A_0[T]$.
So we may write $\partial^b=\sum_{b'} \alpha_{bb'} \frac{\partial}{\partial T_{b'}}$, with $\alpha_{bb'}\in A$.
Hence $\omega^a=\sum_{a'} \beta_{aa'} dT_{a'}$ with $(\beta_{ab})$ the transposed of the inverse of $(\alpha_{ab})$. We get, using Remark \ref{ulle}, that
$$\begin{array}{rcl}
\langle \omega^a_{r}, \partial_s^b \rangle & = & \langle \sum_{a'} \phi_{\cP}^r(\beta_{aa'}) \cP^{-r} \phi_{\cP}^{r*}(dT_{a'}),
\sum_{b'}\phi_{\cP}^s(\alpha_{bb'}) \left(\frac{\partial}{\partial T_{b'}}\right)_s\rangle\\
\  & \  & \  \\
\  & = & \sum_{a'}\sum_{b'}   \phi_{\cP}^r(\beta_{aa'}) \phi_{\cP}^s(\alpha_{bb'}) \langle \cP^{-r} d(\phi_{\cP}^{r*}T_{a'}),
 \cP^s \frac{\partial}{\partial (\phi_{\cP}^s T_{b'})}\rangle\\
 \  & \  & \  \\
 \  & = & \d_{rs}\sum_{a'}\sum_{b'}   \phi_{\cP}^r(\beta_{aa'}) \phi_{\cP}^s(\alpha_{bb'}) \d_{a'b'}\\
 \  & \  & \  \\
 \  & = & \d_{rs}\d_{ab},\end{array}$$
 which ends the proof.
\qed

\begin{corollary}
\label{volll}
Assume $X$ is a smooth affine scheme over $A_0$ possessing a volume form   $\omega$.
Then the form $\bigwedge_{r \leq n} \omega_r$ is a volume form on $\cJ^n_{\cP}(X)$.
\end{corollary}

\medskip


\begin{remark}
The form $\bigwedge_{r \leq n} \omega_r$ is clearly well defined up to a sign.
In fact, if we interpret it in the following way, it is completely well defined. Let
$\bigwedge_{r\leq n}\Omega_X$ denote the module determined by the property that a linear map
$\bigwedge_{r\leq n}\Omega_X\to M$ is the same as an alternating linear map
$\prod_{r\leq n}\Omega_X\to M$,
where $\prod_{r\leq n}\Omega_X$ denotes the set of functions $\{r;r\leq n\}\to\Omega_X$.
Observe that we do not need to choose an ordering on $\{r;r\leq n\}$.
Then define $\bigwedge_{r \leq n} \omega_r$ to be the image of the function
$r\mapsto \omega_r$ under the universal alternating map
$\prod_{r\leq n}\Omega_X\to \bigwedge_{r\leq n}\Omega_X$.
\end{remark}

\medskip

{\it Proof of Corollary \ref{volll}}.
By Corollary \ref{rece} $\cJ^n_{\cP}(X)$ is smooth over $A_0$. Let $N$ be its relative dimension over $A_0$. Then
$$\bigwedge^N \Omega_{\cJ_{\cP}^n(X)}$$
is a locally free sheaf of rank one and the form $\bigwedge_{r \leq n} \omega_r$ is a section of this sheaf. Let $X=\bigcup_i X_i$ be a covering with principal affine open sets such that $\Omega_{X_i}$
is free for each $i$. Then by Proposition \ref{morelli} $\bigwedge_{r \leq n} \omega_r$ is a volume form on each $\cJ^n_{\cP}(X_i)$ and hence on the union $\bigcup_i \cJ_{\cP}^n(X_i)$. But by Proposition \ref{codim2} the complement of this union in $\cJ^n_{\cP}(X)$ has codimension $\geq 2$. This implies that $\bigwedge_{r \leq n} \omega_r$
is a volume form on the whole of $\cJ^n_{\cP}(X)$.
\qed

\begin{remark}
\label{K3}
 Let $S$ be an  abelian scheme or a  K3 surface over $A_0$ (i.e. a smooth projective scheme of relative dimension $2$ with $H^1(S,\cO)=0$ and trivial canonical bundle $\Omega^2_S$) and fix  a volume form $\omega$ on $S$. (So $\omega$ is well defined up to multiplication by a an element of $A_0^{\times}$.) Let $X \subset S$ be an affine open set. Then  and one can consider on each $\cJ^n_{\cP}(X)$
 the volume form $\bigwedge_{r \leq n} \omega_r$. The latter  will be referred to as the  {\it canonical volume form} on $\cJ_{\cP}^n(X)$ and  (once $\omega$ has been fixed) is well defined up to sign.
\end{remark}

\begin{corollary}
Assume the hypotheses of Proposition \ref{morelli}.
For any $f \in A^n=\cO(\cJ^n_{\cP}(X))$ we have the following formula in $\Omega_{A^n}$:
\begin{equation}
\label{kknia}
df=\sum_{1\leq a \leq m}\sum_{0\leq r\leq n} (\partial^a_{r} f) \omega^a_{r}.\end{equation}
\end{corollary}

\begin{remark}
\label{ressp}
By continuity the formula (\ref{kknia}) continues to hold in $\Omega_{(A^n)^{\widehat{p_k}}}$, for any $f \in (A^n)^{\widehat{p_k}}$.
\end{remark}

Finally we will need the following:

\begin{proposition}
\label{invvaa}
Let $Y \ra X$ be an \'{e}tale finite Galois cover of affine smooth schemes over $A_0$ with Galois group $\Gamma$, let $p$ be a prime, and let $n,i \geq 0$ be  integers. Let $\omega$ be an $i$-form on the $p$-jet space $J^n_p(Y)$ which is $\Gamma$-invariant. Then $\omega$ is the pull-back of a unique $i$-form on the $p$-jet space $J^n_p(X)$.
\end{proposition}

\begin{remark}
 Proposition \ref{invvaa}
fails to be true if the $p$-jet spaces $J^n_p(X)=\cJ_{\{p\}}^n(X)^{\widehat{p}}$ and $J^n_p(Y)=\cJ_{\{p\}}^n(Y)^{\widehat{p}}$ are replaced by the arithmetic jet spaces $\cJ_{\{p\}}^n(X)$ and $\cJ_{\{p\}}^n(Y)$ respectively. An example is given by $A_0=\bZ$, $p\neq 2$,  $Y=\Spec \bZ[y,y^{-1}]$,
$X=\Spec \bZ[x,x^{-1}]$, $x \mapsto y^2$ and the $0$-form $\omega:=yy'$, where $y'=\d_p y$. Indeed $yy'$ is invariant under $y \mapsto -y$ but $yy'$ does not belong to the image of
$$\cO(\cJ^1_{\{p\}}(X))\ra \cO(\cJ^1_{\{p\}}(Y)),$$
as one can easily see by tensoring with $\bQ$. On the other hand $yy'$ belongs to the image of
$$ \cO(J^1_{p}(X))\ra \cO(J^1_{p}(Y)),$$
as one can see directly from the formula
$$yy'=x^{(p+1)/2} \sum_{n \geq 2} \left( \begin{array}{c} 1/2 \\ n \end{array} \right) p^{n-1}
\left( \frac{x'}{x^p} \right)^n \in \bZ_p[x,x',x^{-1}]^{\widehat{p}}=\cO(J^1_p(X)).$$

Also  Proposition \ref{invvaa}
fails to be true if the $p$-jet spaces $J^n_p(X)=\cJ_{\{p\}}^n(X)^{\widehat{p}}$ and $J^n_p(Y)=\cJ_{\{p\}}^n(Y)^{\widehat{p}}$ are replaced by  the spaces $\cJ_{\cP}^n(X)^{\widehat{p}}$ and $\cJ_{\cP}^n(Y)^{\widehat{p}}$ respectively, where $\cP$ consists of at least $2$ primes one of which is $p$.
\end{remark}

{\it Proof of Proposition \ref{invvaa}.}
For $i=0$ this was proved in \cite{book}, Proposition 3.27. So we may assume $i \geq 1$.
Let $Y=\bigcup_j Y_j$ be an affine open cover such that $\Omega_{Y_j}$ is free for each $j$, and let $X_j$ be the preimage of $Y_j$ in $X$. Since $J_p^n(X)=\bigcup_j J^n_p(X_j)$ it is sufficient to prove the Proposition for $X_j \ra Y_j$ for each $j$. Replacing $Y$ be $Y_j$ we may assume $\Omega_Y$ is free.
Let $\eta^1,...,\eta^m$ be a basis of $\Omega_Y$ and $\omega^1,...,\omega^m$ be the pull-back of this basis on $X$ which is a basis of $\Omega_X$. By Proposition \ref{morelli} we may write
$$\omega=\sum_{\alpha_1...\alpha_i} a_{\alpha_1...\alpha_i} \omega^{\alpha_1}...\omega^{\alpha_i},$$
for unique  $a_{\alpha_1...\alpha_i} \in \cO(J^n_p(Y))$. Since $\omega$ and $\omega^1,...,\omega^m$ are $\Gamma$-invariant it follows that $a_{\alpha_1...\alpha_i}$ are $\Gamma$-invariant hence (by the $i=0$ case of the Proposition) they are pull-backs of unique functions in $\cO(J^n_p(X))$, which ends our proof.
\qed

\section{The multiplicative group}

In this section we let $A_0=S^{-1}\bZ$ where $S$ is the multiplicative system of all integers coprime to
$p_1,\dots,p_d$. Assume $X=\bG_m=\Spec A^0$, $A^0:=A_0[x,x^{-1}]$, is the multiplicative group scheme
over $A_0$. Then the origin is a uniform point of $X$ with uniform coordinate $T=x-1$ in $X$.
  Let $\omega=\frac{dx}{x} \in \Omega_{\bG_m/A_0}$. Clearly $\omega$ is closed but not exact on
$\bG_m$. Recall that we set $e=(1,\dots,1)\in \bZ_{\geq 0}^d$ and define the $1$-form
	\begin{equation} \label{omega-e-def}
	\omega^{(e)}=-\sum_{0 \leq r \leq e} (-1)^{|r|} \omega_r \in
		H^0(\cJ_{\cP}^e(\bG_m),\Omega_{\cJ^e_{\cP}(\bG_m)}),		
	\end{equation}
where $|r|$ is, as usual, the sum of the
components of $r$ and $\omega_r$ are defined as in (\ref{defome}). Note that $$\omega^{(e)}=-\left(
\prod_{k=1}^d\left( 1-\frac{\phi^*_{p_k}}{p_k}\right) \right)\omega= -\sum_{m|p_1\cdots
p_d}\mu(m)\omega_{[m]},$$ where $\mu$ is the Moebius function and $\omega_{[m]}$ are as in (\ref{sqr}).

Consider now the elements
	$$
	\psi^1_{p_k} \in (A^{e_k})^{\widehat{p_k}}=\bZ_{p_k}[x,x^{-1},\d_{p_k}x]^{\widehat{p_k}},
	$$
	$$\psi^1_{p_k}:=\frac{1}{p_k}
		\log \Big( 1+p \frac{\d_{p_k} x}{x^{p_k}}\Big) =\frac{\d_{p_k}x}{x^{p_k}}-
		\frac{p_k}{2} \Big(\frac{\d_{p_k}x}{x^{p_k}}\Big)^2-\cdots.
	$$
Symbolically, one might write $\psi_{p_k}^1=p_k^{-1}\log \phi_{p_k}(x)- \log x$.
We pass to the multiple-prime case by defining the elements
	\begin{equation} \label{def-f-k}
	f_k:=\bigg(\prod_{\stackrel{l=1}{l\neq k}}^d \Big( 1-\frac{\phi_{p_l}}{p_l}\Big)\bigg)
		\psi^1_{p_k}\in (A^e)^{\widehat{p_k}}.
	\end{equation}
Note that the vector
	\begin{equation} \label{def-psi-e-m}
		\psi^e_m:=(f_1,\dots,f_d)
	\end{equation}
is the {\it arithmetic Laplacian} of $\bG_m$
introduced in \cite{laplace}.
Also consider the logarithm of the formal multiplicative group
	$$
	l_{\bG_m}(T)=\sum_{n=1}^{\infty} (-1)^{n-1}\frac{T^n}{n}\in \bQ[[T]]
	$$
and the series
	\begin{equation} \label{def-psi-e-m-0}
		\psi^e_{m,0}:=-\bigg(\prod_{l=1}^d\Big(1-\frac{\phi_{p_l}}{p_l}\Big)\bigg)
			l_{\bG_m}(T) \in \bQ[[\d_{\cP}^i T; i \leq e]].		
	\end{equation}
(In fact, we have $\psi^e_{m,0}\in A_0[[\d_{\cP}^i T; i \leq e]]$.
See the proof of Theorem 3.3 in \cite{laplace}.)

As discussed in the introduction,
we can think of the $f_k$ and $\psi^e_{m,0}$ as being series expansions
of the formal expression $-\prod_{l=1}^d(1-\frac{\phi_{p_l}}{p_l})\cdot \log(x)$
in different regions of the $A_0$-scheme $\bG_m$. In \cite{laplace}, it was shown that the expansions
agree on what could be interpreted as the intersection these regions. In Theorem \ref{tmu} below,
we show that they are solutions to a common differential equation.

\begin{lemma} \label{lem:dpsi}
$d\psi^1_{p_k} = -\left(1-\frac{\phi^*_{p_k}}{p_k}\right)\omega.$
\end{lemma}

Essentially this was given in \cite{book},
Proposition 7.26, but for the convenience of the reader, we
will include a simple proof here.

\begin{proof}
Let us abbreviate $p=p_k$.  Then we have
\begin{align*}
	p d\psi^1_{p}	&= d\log(1+px^{-p}\delta_{p}x) \\
					&= \frac{d(1+px^{-p}\delta_{p}x)}{1+px^{-p}\delta_{p}x} \\
					&= \frac{d(1+px^{-p}\delta_{p}x)}{1+px^{-p}\delta_{p}x} +
							\frac{d(x^{p})}{x^{p}} - p\frac{dx}{x} \\
					&= \frac{d(x^p+p\delta_{p}x)}{x^p+p\delta_{p}x} - p\frac{dx}{x} \\
					&= \frac{d\phi_{p}(x)}{\phi_{p}(x)} - p\frac{dx}{x} \\
					&= (\phi_{p}^*-p)\frac{dx}{x}.
\end{align*}
Dividing by $p$ then completes the proof.
\end{proof}

Let us say that a $1$-form $\tilde{\omega}$ on some space $\cJ_{\cP}^n(\bG_m)$ is {\it invariant}
if
$$\mu^* \tilde{\omega}=pr_1^* \tilde{\omega}+pr_2^* \tilde{\omega},$$
where
$$\mu,pr_1,pr_2:\cJ_{\cP}^n(\bG_m)\times \cJ_{\cP}^n(\bG_m) \ra \cJ_{\cP}^n(\bG_m)$$
are the multiplication and the $2$ projections respectively. (This agrees with  the standard definition of invariance of $1$-forms because we are in a situation when our group scheme is commutative.)

\begin{lemma}
\label{deas}
The $A_0$-module consisting of the invariant $1$-forms on $\cJ_{\cP}^n(\bG_m)$
is free with basis $\{\omega_r\ |\ r \leq n\}$.
Furthermore, any invariant $1$-form on $\cJ_{\cP}^n(\bG_m)$ is closed.
\end{lemma}

\begin{proof}
	Let us first show each $\omega_r$ is invariant.
	By definition, we have $\omega_r = \varphi_{nr}^*\frac{\phi^{r*}}{\cP^r}(\omega)$.
	And by Proposition \ref{pro:jet-group-maps}, both $\varphi_{nr}$ and $\phi^{r*}$
	are group homomorphisms. Since $\omega$ is an invariant differential on $\bG_m$, we have that
	$\omega_r$ is an invariant differential on $\cJ_{\cP}^n(\bG_m)$.
	
	By Proposition \ref{morelli}, the set $\{\omega_r\ |\ r \leq n\}$ is a
	basis of the $\cO(\cJ_{\cP}^n(\bG_m))$-module of global $1$-forms on $\cJ_{\cP}^n(\bG_m)$.
	In particular, it is $A_0$-linearly independent. It remains to show every invariant
	differential is in its $A_0$-linear span. Let $\tilde{\omega}$ be an invariant differential,
	and write $\tilde{\omega}=\sum_{r \leq n} a_r \omega_r$
	with $a_r \in \cO(\cJ_{\cP}^n(\bG_m))$. Since $\tilde{\omega}$ and the $\omega_r$ are invariant,
	each function $a_r$ is invariant. Therefore $a_r\in A_0$.
	
	Finally, since $\omega$ is closed, so is each $\omega_r$, by (\ref{d-commute}).
	It then follows from the above that every invariant $1$-form is closed.
\end{proof}


Let $P \in \bG_m(A_0)$ be the origin and let $\cS:=\cJ_{\cP}^e(\bG_m)(A_0)$. The canonical lift $P^e \in \cS$
is the origin of $\cJ_{\cP}^e(\bG_m)$. Since $P$ is uniform, $P^e$ is also uniform. Since
$\cJ_{\cP}^e(\bG_m)$ is a group all points in $\cS$ are uniform. Here is our main result on $\bG_m$:

\begin{theorem}
\label{tmu}\

 1) $\omega^{(e)}$ is invariant and hence closed;

 2) $\omega^{(e)}$ is not exact;


 3) The arithmetic Laplacian $\psi^e_m$ (of (\ref{def-psi-e-m}))
 is the $\cP$-adic primitive of $\omega^{(e)}$ normalized along
 $P^e$, and $\psi^e_{m,0}$ (of (\ref{def-psi-e-m-0}))
 is the $P^e$-adic primitive of $\omega^{(e)}$ normalized along $P^e$.
 In particular
 $\omega^{(e)}$ is $\cP$-adically exact and  $\cS$-adically exact.

 4) If $\tilde{\omega}$ is a $1$-form on $\cJ^e_{\cP}(\bG_m)$ that is invariant and $\cP$-adically exact
 then $\tilde{\omega}$ is an $A_0$-multiple of $\omega^{(e)}$.
\end{theorem}

{\it Proof}. Assertion 1) follows from Lemma \ref{deas} and the definition (\ref{omega-e-def})
of $\omega^{(e)}$.

To prove assertion 2) consider the isomorphism
	\begin{equation} \label{isomm}
		\cJ^e_{\cP}(\bG_m) \otimes \bQ\longlabelmap{} \prod_{s \leq e}\bG_{m/\bQ}.
	\end{equation}
induced by $\kappa_{\leq e}$ of (\ref{eq:coghost-map}). Under this map, $\omega_s$ corresponds to the
differential $\cP^{-s}\pr_s^*(\omega)$ on the right-hand side. In particular, if $i$ denotes the
inclusion of $\bG_m$ into the factor $s=0$ (and $1$ on all other factors), then $i^*$ applied to the
differential corresponding to $\omega^{(e)}$ is $\omega$, which is not exact. And so $\omega^{(e)}$
cannot be exact.


To prove the first assertion in 3), we must prove, first, that if $f_k$ is as in
(\ref{def-f-k}), then we have $d f_k=\omega^{(e)}$ for all $k=1,\dots,d$
and, second, that $d \psi^e_{m,0}=\omega^{(e)}$.
By Lemma~\ref{lem:dpsi}, we have
		\begin{align*}
			df_k &=
			d\Bigg(\bigg(\prod_{\stackrel{l=1}{l\neq k}}^d \Big( 1-\frac{\phi^*_{p_l}}{p_l}\Big) \bigg) \psi^1_{p_k}\Bigg) \\
				 &= \bigg(\prod_{\stackrel{l=1}{l\neq k}}^d \Big( 1-\frac{\phi^*_{p_l}}{p_l}\Big)  \bigg) d\psi^1_{p_k} \\
				 &= -\bigg(\prod_{\stackrel{l=1}{l\neq k}}^d \Big( 1-\frac{\phi^*_{p_l}}{p_l}\Big) \bigg) \cdot
					\Big( 1-\frac{\phi^*_{p_k}}{p_k}\Big) \omega\\
				 &= -\bigg(\prod_{l=1}^d \Big( 1-\frac{\phi^*_{p_l}}{p_l}\Big)\bigg) \omega \\
				 &= \omega^{(e)}.
		\end{align*}

The statement about $\psi_{m,0}^e$ follows in the same way. To prove the second assertion in 3)
we use the translation by points in $\cS$ to reduce to the case of the origin.

To prove assertion 4)  embed
	$$
	A_k:=\cO(\cJ_{\cP}^e(\bG_m))^{\widehat{p_k}}
	$$
into
	$$
	B_k:=\bZ_{p_k}[[\d_{\cP}^r T; r\leq e]]
	$$
for each $k$. By hypothesis $\tilde{\omega}=d g_k$, with $g_k \in A_k$ for all $k$. We may assume that
the
$g_k$, viewed as elements of $B_k$, have no constant term. Since $\tilde{\omega}$ is defined over $A_0$
it follows that each $g_k$ belongs to $B_0:=\bQ[[\d_{\cP}^r T; r \leq e]]$ and that the $g_k$ are equal to
some $g_0 \in B_0$. So $g_0\in A_0[[\d_{\cP}^r T;r \leq e]]$. Since
$\mu^* \tilde{\omega}=pr_1^*\tilde{\omega}+pr_2^*\tilde{\omega}$ and $g_0$ has no constant term, it
follows that
	\begin{equation}
		\label{dda}
		\mu^* g_0=pr_1^* g_0+pr_2^* g_0,
	\end{equation}
where
	$$
	\mu^*, pr_1^*, pr_2^*: A_0[[\d_{\cP}^r T;r \leq e]]\ra
	A_0[[\d_{\cP}^r T_1, \d_{\cP}^r T_2;r \leq e]]
	$$
in (\ref{dda}) are induced by the formal group law and the $2$ projections respectively. We conclude that
the tuple $(g_k)$ is a $\d_{\cP}$-character on $\bG_m$ in the sense of Definition 2.33 in \cite{laplace}.
Since the order of $(g_k)$ is $e$, which is also the order of $(f_k)$,
Theorem 3.4 in \cite{laplace} implies there exists $\rho \in A_0$ such
that $g_k=\rho \cdot f_k$. Hence $\tilde{\omega}=\rho \cdot \omega^{(e)}$ and we are done.
\qed

%

\medskip

Consider next, for any multi-index $n$, the volume form $\bigwedge_{r \leq n} \omega_r$ on
$\cJ^n_{\cP}(\bG_m)$. (See Corollary \ref{volll}.)
It will be referred to as the {\it canonical volume form} on $\cJ^n_{\cP}(\bG_m)$.

\begin{corollary}
If $n \geq e$, the canonical volume form on $\cJ_{\cP}^n(\bG_m)$ is $\cP$-adically exact and $\cS$-adically exact but not exact.
\end{corollary}

{\it Proof}. Indeed we have
$\bigwedge_{r \leq n} \omega_r=\omega^{(e)} \wedge \big( \bigwedge_{0 \neq r \leq n} \omega_r \big)$,
which is $\cP$-adically exact and $\cS$-adically exact because $\omega^{(e)}$ is $\cP$-adically exact and $\cS$-adically exact, by Theorem \ref{tmu}.
Assume now that the canonical volume form on $\cJ_{\cP}^n(\bG_m)$ is exact and let us derive a contradiction. Using the isomorphism (\ref{isomm})
we deduce that the form
$$\nu=\frac{dx_1}{x_1} \wedge \cdots \wedge \frac{dx_{N}}{x_{N}}$$
on $\bG_{m/\bQ}^{N}=\Spec \bQ[x_1^{\pm},\dots,x_{N}^{\pm}]$ (where $N$ is the number of elements $r$ that are $\leq n$) is exact. So $\nu=d\eta$,
$$\eta=\sum_{i=1}^{N} f_i \frac{dx_1}{x_1} \wedge \cdots\wedge \widehat{\frac{dx_i}{x_i}} \wedge \cdots\wedge \frac{dx_{N}}{x_{N}},$$
$f_i\in \bQ[x_1^{\pm},\dots,x_{N}^{\pm}]$. Hence
$$\sum_{i=1}^{N} (-1)^i x_i \frac{\partial f_i}{\partial x_i} =1.$$
But this is impossible because none of the Laurent polynomials $x_i\frac{\partial f_i}{\partial x_i}$
has a constant term.
\qed

\begin{remark}
We expect that if $n \not\geq e$ then the canonical  volume form on $\cJ_{\cP}^{n}(\bG_m)$ is not $\cP$-adically exact. In any case this form is not exact (cf. the proof above that applies to any $N$).
\end{remark}

\begin{corollary} \label{cor1}
Consider the derivation $\partial=x \frac{d}{dx}:\cO(\bG_m)\ra \cO(\bG_m)$ and let $\partial_r$ be the
$r$-conjugates of $\partial$. Then
	$$
	\partial_r f_k=(-1)^{|r|}
	$$
for all $k=1,\dots,d$ and  $r \leq e$.
\end{corollary}

{\it Proof}.
By assertion 3) in Theorem \ref{tmu} and Remark \ref{ressp} we have
	$$
	\sum_{r\leq n} (-1)^{|r|}\omega_r = \omega^{(e)} = df_k = \sum_{r\leq n} (\partial_r f_k)\omega_r.
	$$
By Proposition \ref{morelli}, the $\omega_r$ form a basis. Thus $\partial_r f_k = (-1)^{|r|}$.
\qed

\begin{remark}
Using arguments from \cite{laplace} it is easy to show that the $1$-form $\omega^{(e)}$ on $\cJ_{\cP}^e(\bG_m)$ typically has non-zero periods. Indeed consider the
cycle
$$\Gamma=(P^e_1,P^e_2,P^e_3,P^e_4,P^e_1)$$
 on $\cJ_{\cP}^e(\bG_m)$, where $P_j^e \in
  \cJ_{\cP}^e(\bG_m)(A_0)\subset \cJ_{\cP}^e(\bG_m)(\bZ_{p_{k_j}})$
  are the canonical lifts of $P_j \in \bG_m(A_0)$, $k_1=k_2\neq k_3=k_4$,
 $P_2$ and $P_3$ are induced by the identity section  $P_{23}=1 \in \bG_m(A_0)=A_0^{\times}$,
and $P_1$ and $P_4$  are induced by a section $P_{14}\in A_0^{\times}$.   If $P_{14}=\pm 1$, then
	$$
	\int_{\Gamma} \omega^{(e)} =0.
	$$
Indeed, each of the  $p_{k_j}$-adic primitives $f_{k_j}$ for $\omega^{(e)}$ gives a group homomorphism $\cJ_{\cP}^e(\bG_m)(\bZ_{p_{k_j}})\to \bG_a(\bZ_{p_{k_j}})$.
But the target is torsion free, and $P_{14}^e$ is a torsion point. So we have $f_{k_j}(P_{14}^e)=0$.

However, if $P_{14} \neq \pm 1$,  we have that
	$$
	\int_{\Gamma} \omega^{(e)} = f_{k_1}(P^e_{14})-f_{k_4}(P^e_{14}) \neq 0.
	$$
This follows from  \cite{laplace}, proof of Theorem 3.4.
(The argument is that
$f_{k_1}(P^e_{14})$ is a non-zero rational number times the $p_{k_1}$-adic logarithm of an element
in $1+p_{k_1}\bZ_{p_{k_1}} \setminus \{1\}$; but the latter logarithm is not in $\bQ$ by
Mahler's $p$-adic analogue \cite{mahler, bertrand} of the Hermite-Lindeman theorem.)
\end{remark}

\section{Elliptic curves}
Again we let $A_0=S^{-1}\bZ$ where $S$ is the multiplicative system of all integers coprime to
$p_1,\dots,p_d$. Assume all primes in $\cP$ are $\geq 5$.
Consider an elliptic curve over $A_0$,
	$$
	E:=E_{a,b}:=\Proj A_0[x_0,x_1,x_2]/(x_0x_2^2-x_1^3-ax_1x_0^2-bx_0^3),
	$$
with $a,b\in A_0$ and
$-4a^3-27b^2 \in A_0^{\times}$. Let $\omega=\frac{dx}{y}\in H^0(E,\Omega_E)$,
$x:=\frac{x_1}{x_0}$, $y:=\frac{x_2}{x_0}$, $T=\frac{x}{2y}$. Let $X \subset E$ be the affine open set
where 
$x_2$ is invertible.
The origin $P=[0,0,1]$ is uniform with uniform coordinate $T$.
We continue to denote
by $\omega$ the image of $\omega$ in $H^0(X,\Omega_X)$. Clearly $\omega$ is closed but not exact on $X$;
indeed any non-zero exact form on $X$ must have a pole in $E \setminus X$. Let $a_{p_k}$ be the trace of
Frobenius on $E \otimes {\mathbb F}_{p_k}$. (So $E$ has $1-a_{p_k}+p_k$ points with
coordinates in ${\mathbb F}_{p_k}$.)
Also we extend this definition by setting
$$a_m:=a_{p_1}^{i_1}\cdots a_{p_d}^{i_d}$$
for $m=p_1^{i_1}\cdots p_d^{i_d}|p_1\cdots p_d$, $i_j\in \{0,1\}$.
Consider the form
	$$
	\omega^{(2e)}:= \left( \prod_{l=1}^d \bigg( 1-a_{p_l}
	\frac{\phi^*_{p_l}}{p_l}+p_l\frac{\phi_{p_l^2}^*}{p_l^2} \bigg) \right) \omega
	\in H^0\big(\cJ_{\cP}^{2e}(X),\Omega_{\cJ_{\cP}^{2e}(X)}\big).
	$$
One can write
	$$
	\omega^{(2e)}=\sum_{m|p_1^2\cdots p_d^2} \mu(m')m'' a_{m'}  \omega_{[m]},
	$$
where $m=m'(m'')^2$ with $m', m''$ square free and coprime;
$\mu$ is the Moebius function; and $\omega_{[m]}$ is as in (\ref{sqr}).
It follows  from \cite{book}, Theorem 7.22 and Corollary 7.28, plus \cite{frob} Theorem 1.10, that there
exist elements
$$\psi_{p_k}^2 \in \cO(J_{p_k}^{2}(E))$$
(where $J_{p_k}^2(E)$ denotes the formal scheme defined in Remark \ref{nonaffinejetspaces})
such that $\psi_{p_k}^2$ vanish at $0$ and
\begin{equation}
\label{right}
d\psi^2_{p_k}=
\left(1-a_{p_k}\frac{\phi^*_{p_k}}{p_k}+p_k\left(\frac{\phi_{p_k}^*}{p_k}\right)^2\right)\omega.
\end{equation}
Clearly $\psi^2_{p_k}$ are unique with the above properties.
We continue to denote by $\psi_{p_k}^2$ the image of this element in $\cO(J^2_{p_k}(X))=\cO(\cJ_{\cP}^{2e_k}(X)^{\widehat{p_k}})$. (N.B. The superscript $2$ in $\psi^2_{p_k}$ is not an exponent; it merely indicates that the order of that element is $2$.)
 Following \cite{laplace} consider the elements
$$f_k:=\left( \prod_{\stackrel{l=1}{l \neq k}}^d \left( 1-a_{p_l} \frac{\phi_{p_l}}{p_l}+p_l\frac{\phi_{p_l^2}}{p_l^2} \right) \right) \psi^2_{p_k}\in
\cO(\cJ_{\cP}^{2e}(X)^{\widehat{p_k}}).$$
The vector
\begin{equation}
\label{ale}
\psi^{2e}_E:=(f_1,\dots,f_d)
 \end{equation}
 is the (restriction to $X$ of the) {\it arithmetic Laplacian} of $E$ introduced in \cite{laplace}. Moreover we may consider
 \begin{equation}
 \label{ale0}
\psi^{2e}_{E,0}:=\left[ \prod_{l=1}^d\left(1-a_{p_l}
\frac{\phi_{p_l}}{p_l}+p_l\left(
\frac{\phi_{p_l}}{p_l}\right)^2\right)\right]l_E(T)\in \bQ[[\d_{\cP}^i T;i \leq 2e]],\end{equation}
where $l_E(T) \in \bQ[[T]]$ is the logarithm of the formal group of $E$. By \cite{laplace}, proof of Theorem 3.7, we have
$$\psi_{E,0}^{2e}\in A_0[[\d_{\cP}^i T; i \leq 2e]].$$

Let $\mu:E \times E \ra E$ be the multiplication map and let
$$X^{(2)}:=(X \times X)\cap \mu^{-1}(X) \subset E \times E.$$
Let us say that a $1$-form $\tilde{\omega}$ on some space $\cJ_{\cP}^n(X)$ is {\it invariant} if
$$\mu^* \tilde{\omega}=\pi_1^* \tilde{\omega}+\pi_2^* \tilde{\omega},$$
where
$\mu,\pi_1,\pi_2:\cJ_{\cP}^n(X^{(2)})\ra \cJ_{\cP}^n(X)$
come from the multiplication and the $2$ projections.
(Again this agrees with  the standard definition of invariance of $1$-forms because we are in a commutative situation.)

 \begin{lemma}
 \label{deasu}
The $A_0$-module consisting of the invariant $1$-forms on $\cJ_{\cP}^n(X)$
is free with basis $\{\omega_r\ |\ r \leq n\}$.
Furthermore, any invariant $1$-form on $\cJ_{\cP}^n(X)$ is closed.
\end{lemma}

\begin{proof}
	Exactly as in the case of $\bG_m$.
\end{proof}

Since the origin $P=[0,0,1]\in X(A_0)$ is uniform so is its canonical lift $P^{2e} \in \cJ_{\cP}^{2e}(X)(A_0)$; cf. section 3. Let $\cS=\cJ_{\cP}^{2e}(X)(A_0)$. Since $\cJ_{\cP}^{2e}(X)$ is a group all points in $\cS$ are uniform.

\begin{theorem}
\label{tme}\

 1) $\omega^{(2e)}$ is invariant and hence closed;

 2) $\omega^{(2e)}$ is not exact;

 3) The arithmetic Laplacian $\psi^{2e}_E$  (of (\ref{ale})) is the $\cP$-adic primitive of $\omega^{(2e)}$ normalized along $P^{2e}$ and $\psi^{2e}_{E,0}$ (of (\ref{ale0})) is the $P^{2e}$-adic primitive of $\omega^{(2e)}$ normalized along $P^{2e}$. In particular
 $\omega^{(2e)}$ is $\cP$-adically exact and $\cS$-adically exact.

 4) Assume $E$ has ordinary reduction at all the primes in $\cP$ and  $\tilde{\omega}$ is a $1$-form on $\cJ^{2e}_{\cP}(X)$ that is invariant and $\cP$-adically exact. Then
  $\tilde{\omega}$ is an $A_0$-multiple of $\omega^{(2e)}$.
\end{theorem}

{\it Proof}. Assertions 1) and 2)  follow exactly as in Theorem \ref{tmu}.

To prove assertion 3) we must prove, as in the case of $\bG_m$,  that $d f_k=\omega^{(2e)}$ for all $k=1,\dots ,d$ and $d \psi_{E,0}^{2e}=\omega^{(2e)}$.
By (\ref{right}) we have
		\begin{align*}
			df_k &=
			d\left(\left(\prod_{\stackrel{l=1}{l\neq k}}^d \left( 1-a_{p_l}\frac{\phi_{p_l}}{p_l}+p_l
\left( \frac{\phi_{p_l}}{p_l}\right)^2\right) \right) \psi^1_{p_k}\right) \\
				 &= \left(\prod_{\stackrel{l=1}{l\neq k}}^d \left( 1-a_{p_l}\frac{\phi^*_{p_l}}{p_l}+p_l\left(\frac{\phi^*_{p_l}}{p_l}\right)^2\right)  \right) d\psi^1_{p_k} \\
				 &= \left(\prod_{\stackrel{l=1}{l\neq k}}^d \left( 1-a_{p_l}\frac{\phi^*_{p_l}}{p_l}+p_l\left(\frac{\phi^*_{p_l}}{p_l}\right)^2\right) \cdot
					\left( 1-a_{p_k}\frac{\phi^*_{p_k}}{p_k}+p_k\left(\frac{\phi^*_{p_k}}{p_k}\right)^2\right) \right) \omega\\
				 &= \left(\prod_{l=1}^d \left( 1-a_{p_l}\frac{\phi^*_{p_l}}{p_l}+p_l\left(\frac{\phi^*_{p_l}}{p_l}\right)^2\right)\right) \omega \\
				 &= \omega^{(2e)}.
		\end{align*}

A similar computation proves the statement about $\psi^{2e}_{E,0}$.

To prove assertion 4) we proceed exactly as in the proof of assertion 4) in Theorem \ref{tmu}; instead of Theorem 3.4 in \cite{laplace} we need to use Theorem 3.8 in \cite{laplace} in conjunction with Lemma 7.33 in \cite{book}.
\qed

%

\medskip

Consider next, for any multi-index $n$, the volume form $\bigwedge_{r \leq n} \omega_r$ on
$\cJ^n_{\cP}(X)$. (See Corollary \ref{volll}.)
It will be referred to as the {\it canonical volume form} on $\cJ^n_{\cP}(\bG_m)$.
Also let $\cS=\cJ_{\cP}^n(X)(A_0)$.


\begin{corollary}
If $n \geq 2e$ the canonical volume form on $\cJ_{\cP}^{n}(X)$ is $\cP$-adically exact and $\cS$-adically exact but not exact.
\end{corollary}

{\it Proof}.
The same argument as for $\bG_m$. To prove  non-exactness  it is enough to prove, as in the case of $\bG_m$, that if we view $E$ and $X$ as schemes over $\bQ$ and
$\nu \in H^0(E^N,\Omega^N_{E^N}) \subset H^N_{DR}(E^N)$ is a volume form (where $N\geq 1$ is an integer)
then the image of $\nu$ via the restriction map
$$H^N_{DR}(E^N)\ra H^N_{DR}(X^N)$$
is non-zero. But $\nu$ lies in the $H^1_{DR}(E) \otimes \cdots \otimes H^1_{DR}(E)$ K\"{u}nneth summand
of $H^N_{DR}(E^N)$ so the image of $\nu$ in $H^N_{DR}(X^N)$ lies in the
$H^1_{DR}(X) \otimes \cdots \otimes H^1_{DR}(X)$ K\"{u}nneth summand of $H^N_{DR}(X^N)$. So we are
reduced to show that the restriction map $H^1_{DR}(E)\ra H^1_{DR}(X)$ is injective, which is true.
\qed

\begin{remark}
We expect that if $n \not\geq 2e$ then no invariant volume form on $\cJ_{\cP}^{n}(X)$ is  $\cP$-adically exact. In any case no such form is exact (cf. the proof above that applies to any $N$).
\end{remark}

\begin{remark}
Let $S$ be an abelian scheme or a  K3 surface over $A_0$ and let $X \subset S$ be an affine open set. It is natural to ask if the canonical volume forms (cf. Remark \ref{K3}) on the various jet spaces $\cJ^n_{\cP}(X)$
are $\cP$-adically exact. We are not able to treat this question in general. Note however that by Theorem \ref{tme} the answer to this question is positive for abelian schemes of the form  $S=E_1 \times \cdots  \times E_m$ with $E_j$  elliptic curves as in this section and $X=X_1 \times \cdots  \times X_m$ with $X_j \subset E_j$ affine open sets. A special case of K3 surfaces will be treated in the next section.
\end{remark}

\begin{corollary}
 Consider the derivation $\partial=y \frac{\partial}{\partial x}:\cO(X)\ra \cO(X)$ and let $\partial_r$ be the $r$-conjugates of $\partial$. Then $$\partial_r f_k=\mu(m')m'' a_{m'}$$
 for all $k=1,\dots ,d$ and  $r \leq 2e$, $m=\cP^r$.
\end{corollary}

{\it Proof}.
As in the Corollary \ref{cor1}, this follows from assertion 3) in Theorem \ref{tme} and Remark
\ref{ressp}.
\qed

\begin{remark}
As in the case of $\bG_m$, using arguments from \cite{laplace}, it is easy to show that the $1$-form
$\omega^{(2e)}$ on $\cJ_{\cP}^{2e}(X)$ typically has non-zero periods. Indeed consider the cycle
$$\Gamma=(P^{2e}_1,P^{2e}_2,P^{2e}_3,P^{2e}_4,P^{2e}_1)$$
on $\cJ_{\cP}^{2e}(X)$, where $P_j^{2e} \in \cJ_{\cP}^{2e}(X)(A_0) \subset \cJ_{\cP}^{2e}(X)(\bZ_{p_{k_j}})$ are the canonical lifts of $P_j \in X(A_0)$, $k_1=k_2 \neq k_3=k_4$,
$P_2$ and $P_3$  are induced by the identity section  of $E(A_0)$
and $P_1$ and $P_4$ are induced by a section $P_{14} \in X(A_0)$.  As with $\bG_m$,
if $P_{14}$ is torsion then
$$\int_{\Gamma} \omega^{(2e)} =0.$$
However, if $P_{14}$ is non-torsion,  we have that
$$\int_{\Gamma} \omega^{(2e)} = f_{k_1}(P^{2e}_{14})-f_{k_4}(P^{2e}_{14}) \neq 0.$$
This follows from the proof of Theorem 3.8 in \cite{laplace}. (The argument is that
$f_{k_1}(P^{2e}_{14})$ is a non-zero rational number times the $p_{k_1}$-adic elliptic logarithm of an element
in $p_{k_1}\bZ_{p_{k_1}} \setminus \{0\}$; but the latter logarithm is not in $\bQ$ by
Bertrand's $p$-adic analogue \cite{bertrand} of the Hermite-Lindeman theorem.)
\end{remark}

\section{K3 surfaces}

In this section we consider our theory for a special class of K3 surfaces, namely for  Kummer surfaces attached to products of two elliptic curves over $\bQ$. For simplicity we will only discuss $\cP$-adic exactness; an analysis  of $\cS$-adic exactness can be done but  requires a slight generalization of the discussion in the previous section (in which sections are replaced by multisections) so it will be omitted here.

Again we let $A_0=S^{-1}\bZ$ where $S$ is the multiplicative system of all integers coprime to $p_1,\dots ,p_d$ and assume that all primes in $\cP$ are $\geq 5$. Consider two elliptic curves $E=E_{a,b}$ and $\tilde{E}=\tilde{E}_{\tilde{a},\tilde{b}}$
over $A_0$,
as in the previous section, defined by cubics with coefficients $a,b \in A_0$ and $\tilde{a},\tilde{b}\in A_0$ respectively. Denote by $\omega$ and $\tilde{\omega}$ the corresponding $1$-forms on $E$ and $\tilde{E}$ and denote by $a_{p_k}$ and $\tilde{a}_{p_k}$ the corresponding traces of Frobenius.
 For each $k=1,\dots ,d$ we have at our disposal the functions $\psi_{p_k}^2 \in \cO(J^2_{p_k}(E))$ and $\tilde{\psi}_{p_k}^2 \in \cO(J^2_{p_k}(\tilde{E}))$
respectively. Consider the abelian scheme $A=E \times \tilde{E}$. We view $\omega$ and $\tilde{\omega}$ as $1$-forms on $A$ via pull-back. Then $\alpha:=\omega \wedge \tilde{\omega}$ is a volume form on $A$.  We also view $\psi^2_{p_k}$ and $\tilde{\psi}^2_{p_k}$ as elements
of $\cO(J^2_p(A))$ via pull-back.

Let $T \subset E$ and $\tilde{T} \subset \tilde{E}$ be the kernels of $[2]$ (the multiplication by $2$)
on $E$ and $\tilde{E}$ respectively and assume that $T$ and $\tilde{T}$ are unions of sections of our
elliptic curves over $A_0$ (in other words we assume that the cubics $x^3+ax+b$ and
$\tilde{x}^3+\tilde{a}\tilde{x}+\tilde{b}$ have all their roots in $A_0$). Let $S=Km(A)$ be the Kummer
surface attached to $A$, i.e., $S=B/i$ where $B$ is the blow up of $A$ at
$T \times \tilde{T}$ and $i:B \ra B$ is the involution lifting of the multiplication by $-1$,
$[-1]:A\ra A$. We have a diagram
$$A \stackrel{\epsilon}{\longleftarrow} B \stackrel{\pi}{\longrightarrow} S.$$
Recall that
$$\pi(\epsilon^{-1}(T \times \tilde{T}))=\bigcup_{j=1}^{16} L_j$$
where $L_j$ are  $\bP^1$s on $S$ with self-intersection $-2$.
Let $\sigma$ be the volume form on $S$
normalized such that $\pi^* \sigma=\epsilon^* \alpha$ on $B$.

Let $U \subset A$ be an affine open set that is {\it symmetric} (i.e. $[-1]U=U$) and disjoint from $T \times \tilde{T}$. Moreover let $Y=\epsilon^{-1}(U)$ and $X=Y/i$.
So we have a diagram
$$U \stackrel{\epsilon}{\longleftarrow} Y \stackrel{\pi}{\longrightarrow} X$$
where the first map is an isomorphism and the second is a finite \'{e}tale Galois cover with group $\langle i \rangle$.

For any multi-index $n\in \bZ_{\geq 0}^d$ consider the {\it canonical volume form}
$\bigwedge_{r \leq n} \sigma_r$ on $\cJ_{\cP}^n(X)$.

\begin{theorem}
\label{K3thm}
Let $n \geq 2e$. Then the canonical volume form  on $\cJ_{\cP}^n(X)$ is $\cP$-adically exact. Moreover if $U=(E\setminus T)\times (\tilde{E}\setminus \tilde{T})$, then the
canonical volume form  on $\cJ_{\cP}^n(X)$ is not exact.
\end{theorem}

\medskip

In order to prove our Theorem we need to introduce certain $1$-forms.
 Indeed, for each $k=1,\dots ,d$ consider the $1$-form
\begin{equation}
\label{etas}
\eta_k  :=   \psi^2_{p_k} \tilde{\omega}_0 = \psi^2_{p_k}\tilde{\omega}
\end{equation}
on $J_{p_k}^2(A)$ and hence on $J_{p_k}^2(U)$. Note that
\begin{equation}
\label{exppr}
d \eta_k  = d \psi_{p_k}^2 \wedge \tilde{\omega}_0=\omega_0 \wedge \tilde{\omega}_0-a_{p_k} \omega_{e_k} \wedge \tilde{\omega}_0+p_k \omega_{2e_k} \wedge \tilde{\omega}_0.
\end{equation}
Clearly $[-1]^* \psi_{p_k}^2=-\psi^2_{p_k}$ and $[-1]^*\tilde{\omega}_0=-\tilde{\omega}_0$ on $J_{p_k}^2(A)$ so $[-1]^* \eta_k=\eta_k$ and hence $i^* (\epsilon^* \eta_k)=\epsilon^* \eta_k$ on $J^2_{p_k}(Y)$. By Proposition \ref{invvaa}
\begin{equation}
\label{theta}
\epsilon^* \eta_k= \pi^* \theta_k
\end{equation}
on $J^2_{p_k}(Y)$
for some $1$-form $\theta_k$ on $J_{p_k}^2(X)$.

\medskip

{\it Proof of Theorem \ref{K3thm}.}
Let us still denote by $\theta_k$ the induced $1$-form on $\cJ_{\cP}^n(X)^{\widehat{p_k}}$. We claim that
\begin{equation}
\label{eqq}
\bigwedge_{r \leq n} \sigma_r=d \theta_k \wedge  \bigwedge_{0 \neq r \leq n} \sigma_r =
d \Big( \theta_k \wedge  \bigwedge_{0 \neq r \leq n} \sigma_r  \Big)
\end{equation}
on $\cJ_{\cP}^n(X)^{\widehat{p_k}}$. This will prove that $\bigwedge_{r \leq n}\sigma_r$ is
$\cP$-adically exact on $\cJ_{\cP}^n(X)$. The second equality in (\ref{eqq}) is clear because $\sigma$ is
closed and hence each $\sigma_r$ is closed. To check the first equality in (\ref{eqq}), it is enough to
prove that the left-hand side and the right-hand side become equal when pulled back by $\pi$,
since $\pi$ is \'etale. Now on the one hand, we have
$$
\pi^* \Big( \bigwedge_{r \leq n} \sigma_r\Big) = \bigwedge_{r \leq n} \pi^* \sigma_r
 	= \bigwedge_{r \leq n} (\pi^* \sigma)_r
	= \bigwedge_{r \leq n} (\epsilon^* \alpha)_r
 	= \epsilon^* \Big( \bigwedge_{r \leq n} \alpha_r\Big).$$
Similarly
\begin{equation}
\label{similarly}
\begin{array}{c}
\pi^* \left( \bigwedge_{0 \neq r \leq n} \sigma_r\right)=\epsilon^* \left( \bigwedge_{0 \neq r \leq n} \alpha_r\right).\end{array}\end{equation}
On the other hand we have
$$\begin{array}{rcll}
\pi^* \left( d  \theta_k \wedge \big( \bigwedge_{0 \neq r \leq n} \sigma_r \big) \right) & = &
d (\pi^* \theta_k) \wedge \pi^* \left( \bigwedge_{0 \neq r \leq n} \sigma_r \right) & \\
\  & \  & \ & \ \\
\  & = & d(\epsilon^* \eta_k) \wedge \epsilon^*\left( \bigwedge_{0 \neq r \leq n} \alpha_r\right) & \text{by (\ref{similarly})}\\
\ & \ & \  & \ \\
\ & = & \epsilon^* \left( d \eta_k \wedge \big( \bigwedge_{0 \neq r \leq n} (\omega_r \wedge \tilde{\omega}_r) \big) \right) & \ \\
\  & \  & \ & \   \\
\  & = & \epsilon^* \left( \omega_0 \wedge \tilde{\omega}_0 \wedge \big( \bigwedge_{0 \neq r \leq n}(\omega_r \wedge \tilde{\omega}_r) \big) \right) & \text{by (\ref{exppr})}\\
\  & \  & \ & \   \\
\  & = & \epsilon^* \left( \bigwedge_{r \leq n} \alpha_r\right), & \
\end{array}$$
which ends the proof of the $\cP$-adic exactness assertion in the Theorem. To prove the second
assertion of the Theorem,
let us assume for a contradiction that $\bigwedge_{r \leq n}\sigma_r$ is exact on $\cJ_{\cP}^n(X)$.
Pulling back by $\pi$, we get that $\bigwedge_{r \leq n}\alpha_r$ is exact on $\cJ_{\cP}^n(U)$. Tensoring
with $\bQ$, we get that the volume form on $E^N \times \tilde{E}^N$ over $\bQ$ restricted to
$(E \setminus T)^N \times (\tilde{E}\setminus \tilde{T})^N$ is exact for some $N$. But this cannot be
the case, by the same argument as in the proof of Theorem \ref{tme}.
\qed

\bigskip

Of course the affine scheme $X \subset S$ we have been considering is disjoint from the union $\bigcup_{j=1}^{16} L_j$ of the $16$ lines on $S$. On the other hand we expect that if $X' \subset S$ is an affine open set which intersects $\bigcup_{j=1}^{16} L_j$,
then the volume form $\bigwedge_{r \leq 2e} \sigma_r$
on $\cJ_{\cP}^{2e}(X')$ is not $\cP$-adically exact. We can only prove a partial result in this direction; see Theorem \ref{sixforms} below.

By the way each $1$-form $\theta_k$ extends to a $1$-form on
$$J^2_{p_k}\left(S \setminus \bigcup_{j=1}^{16}L_j\right)$$ simply because $S \setminus
\bigcup_{j=1}^{16}L_j$ can be covered by open subsets $X$ of the type we have been considering. Now if
the $1$-forms $\theta_k$ could be extended to $1$-forms on the whole of $J^2_{p_k}(S)$ then it would
follow (as in the proof of Theorem \ref{K3thm}) that the canonical volume form on $\cJ_{\cP}^{2e}(X')$ is
$\cP$-adically closed for any affine open set $X' \subset S$. As mentioned above we do not expect this to
be true. And indeed the extension property for $\theta_k$ does not hold, as shown by the following
Theorem. Let us fix $k$ in $\{1,\dots ,d\}$

\begin{theorem}
\label{sixforms}
 The  $1$-form $\theta_k$ on $J^2_{p_k}\left(S \setminus \bigcup_{j=1}^{16}L_j\right)$  cannot be  extended to a $1$-form on the whole of $J^2_{p_k}(S)$.
\end{theorem}

{\it Proof.}
Let us write
write $p=p_k$, $\psi_p^2=\psi_{p_k}^2$, $\eta=\eta_k$, $\theta=\theta_k$ (cf. (\ref{theta})), etc.
 Assume $\theta$ can be  extended to a $1$-form on $J^2_{p}(S)$.
 Let $z=x/y$ $\tilde{z}=\tilde{x}/\tilde{y}$. Then the completion of $A$ along the origin has ring (of global sections) $A_0[[z,\tilde{z}]]$.  There is a point $P$ on $B$,
lying above the origin of $A$ such that the completion of $B$ along $P$ has ring $A_0[[z,v]]$ where $\tilde{z}=zv$. Then the completion of $S$ along the image of $P$ has ring $A_0[[u,v]]$, with $u=z^2$.
Note that there is a natural identification
$$\bQ_p[[z,\tilde{z},z',\tilde{z}', z'', \tilde{z}'']]=\bQ_p[[z,\tilde{z},z^{\phi}, \tilde{z}^{\phi},z^{\phi^2},\tilde{z}^{\phi^2}]],$$
where $z'=\d_{p}z, z''=\d_{p}^2 z$, etc., and $z^{\phi}=\phi_{p}(z)=z^{p}+p z'$, $z^{\phi^2}=\phi_{p}^2(z)$, etc. So we have natural inclusions
$$\bQ_p[[z,\tilde{z},z^{\phi}, \tilde{z}^{\phi},z^{\phi^2},\tilde{z}^{\phi^2}]] \rightarrow
\bQ_p[[z,v,z^{\phi}, v^{\phi},z^{\phi^2},v^{\phi^2}]] \leftarrow
\bQ_p[[u,v,u^{\phi}, v^{\phi},u^{\phi^2},v^{\phi^2}]].$$
Let $\sum_{n \geq 1} c_n z^n$ be the logarithm of the formal group of $E$ with respect to $z$ and let
$\sum_{n \geq 1} \tilde{c}_n \tilde{z}^n$ be defined similarly. (So $c_1=\tilde{c}_1=1$.) Then
$$\psi_p^2=\frac{1}{p}(\phi^2-a_p\phi+p)(\sum c_nz^n)$$
in $\bQ_p[[z,z^{\phi},z^{\phi^2}]]$ and similarly for $\tilde{\psi}_p^2$. So we have
\begin{equation}
\label{ecu1}
\eta  =  \left(\sum c_nz^n-\frac{a_p}{p} \sum c_n (z^{\phi})^n+\frac{1}{p} \sum c_n (z^{\phi^2})^n\right) \times \left(
\sum n \tilde{c}_n \tilde{z}^{n-1} d \tilde{z}\right).\end{equation}
Since we assumed that $\theta$ extends to $J^2_p(S)$ it follows that there exist series
$$f_0,f_1,f_2, g_0, g_1, g_2 \in \bQ_p[[u,v,u^{\phi}, v^{\phi},u^{\phi^2},v^{\phi^2}]]$$
such that
\begin{equation}
\label{ecu2}
\eta=f_0 du+f_1 d u^{\phi}+f_2du^{\phi^2}+g_0dv+g_1dv^{\phi}+g_2dv^{\phi^2}.\end{equation}
Replacing $\tilde{z}$ by $zv$ in (\ref{ecu1}) and $u$ by $z^2$ in (\ref{ecu2}), expressing $\eta$ as a linear combination of
$$dz, dz^{\phi}, dz^{\phi^2}, dv, dv^{\phi}, dv^{\phi^2},$$
 and picking out the coefficient of $dz$ we get
$$\left(\sum c_nz^n-\frac{a_p}{p} \sum c_n (z^{\phi})^n+\frac{1}{p} \sum c_n (z^{\phi^2})^n\right)
\times \left(\sum n \tilde{c}_n z^{n-1}v^n\right)=f_0 \times 2z$$
of series in $\bQ_p[[z,v,z^{\phi},v^{\phi},z^{\phi^2},v^{\phi^2}]]$.
Picking out the coefficient of $z^{\phi^2}v$ we get $\frac{1}{p}=0$, a contradiction.
\qed

\begin{remark}
 Note that
the form $\theta_k$  is just one of six $1$-forms
\begin{equation}
\label{bunchoforms}
\theta_k, \theta'_k,\theta''_k,\tilde{\theta}_k, \tilde{\theta}'_k, \tilde{\theta}''_k\end{equation}
on $J_{p_k}^2(X)$ that one could consider and that could be used in the proof of Theorem \ref{K3thm}
in the same way in which $\theta_k$ was used. These six $1$-forms are the unique $1$-forms on $J^2_{p_k}(X)$
such that
	\begin{equation} \label{def-thetas}
	\begin{array}{lll}
	\pi^* \theta_k= \epsilon^*(\psi^2_{p_k} \tilde{\omega}_0), &
	\pi^* \theta'_k= \epsilon^*(\psi^2_{p_k} \tilde{\omega}_{e_k}), &
	\pi^* \theta''_k= \epsilon^*(\psi^2_{p_k} \tilde{\omega}_{2e_k}),\\
	\  & \  & \  \\
	\pi^* \tilde{\theta_k}= \epsilon^*(\tilde{\psi}^2_{p_k} \omega_0), &
	\pi^* \tilde{\theta}'_k= \epsilon^*(\tilde{\psi}^2_{p_k} \omega_{e_k}), &
	\pi^* \tilde{\theta}''_k= \epsilon^*(\tilde{\psi}^2_{p_k} \omega_{2e_k}).\end{array}			
	\end{equation}
They can all be extended to $J^2_{p_k}(S \setminus \bigcup_{j=1}^{16}L_j)$. On the other hand, as before, none of them can be extended to $J^2_{p_k}(S)$.

It would be interesting to know if the spaces $\cJ_{\cP}^n(X)$ (where $X$ is as above) possess non-zero $2$-forms that are $\cP$-adically
exact. For instance one can ask the following question: {\it is there a non-trivial $A_0$-linear combination of the $\sigma_r$s (with $r \leq 2e$) on $\cJ_{\cP}^{2e}(X)$ which is $\cP$-adically exact?} The answer to  the latter question is {\it yes} in case $E$ and $\tilde{E}$ have supersingular reduction at all the primes in $\cP$; see Theorem \ref{ultima} below.
But we expect that the answer to this question is {\it no} in case
$E$ and $\tilde{E}$ have ordinary reduction at the primes in $\cP$. Again, see the comments below.

Indeed it is a trivial exercise in linear algebra to show that if $E$ and $\tilde{E}$ have ordinary reduction at $p_k$ (hence $a_{p_k}\neq 0$, $\tilde{a}_{p_k} \neq 0$) then no non-zero $2$-form on $J^{2e_k}_{\cP}(X)$ that is a  $A_0$-linear combination of $\sigma_0, \sigma_{e_k}, \sigma_{e_k}$ can be a $\bZ_{p_k}$-linear combination
of the forms
$$d\theta_k, d\theta'_k,d\theta''_k,d\tilde{\theta}_k, d\tilde{\theta}'_k, d\tilde{\theta}''_k$$
on $J^2_{p_k}(X)$.


On the other hand, if $E$  and $\tilde{E}$ have supersingular reduction
at $p_k$ (hence $a_{p_k}=\tilde{a}_{p_k}=0$), then
\begin{equation}
\label{memo}
\sigma_0-p_k^2\sigma_{2e_k}
=\frac{1}{2}(d\theta_k-p_k d \theta_k''-d \tilde{\theta}_k +p_k d \tilde{\theta}''_k).
\end{equation}
This can be checked by a straight-forward computation using (\ref{right}) and the
defining formulas (\ref{def-thetas}).
\end{remark}

In the next Theorem we consider the $2$-form
$$\sigma^{(2e)}:=\bigg(\prod_{k=1}^d\Big(1-\frac{\phi_{p_k}^{*2}}{p_k^2}\Big)\bigg)\sigma=
\sum_{m|p_1\cdots p_d}\mu(m)m^2\sigma_{[m]}$$
on $\cJ_{\cP}^{2e}(X)$,
where $\mu$ is the Moebius function and $\sigma_{[m]}$ are defined as in (\ref{sqr}).

\begin{theorem}
\label{ultima}
Assume $E$ and $\tilde{E}$ have supersingular reduction at all primes in $\cP$.
Then the $2$-form $\sigma^{(2e)}$
on $\cJ_{\cP}^{2e}(X)$ is $\cP$-adically exact.
Moreover if $U=(E\setminus T)\times (\tilde{E}\setminus \tilde{T})$
 this $2$-form is not exact.
\end{theorem}

{\it Proof}.
Exactly as in the proof of Theorems \ref{tmu} and \ref{tme}, using (\ref{memo}), one shows that
$\sigma^{(2e)}=d\beta_l$
for each $l=1,\dots ,d$, where $\beta_l$ are the $1$-forms
$$\beta_l:=\frac{1}{2} \left(\prod_{\stackrel{k=1}{k\neq l}}^d \left(1-\frac{\phi_{p_k}^{*2}}{p_k^2}\right)\right) (\theta_k-p_k  \theta_k''- \tilde{\theta}_k +p_k  \tilde{\theta}''_k)$$
on $\cJ_{\cP}^{2e}(X)^{\widehat{p_l}}$ and $\cP$-adic exactness follows. Non-exactness follows as in Theorem \ref{K3thm}.
\qed

\bibliographystyle{amsplain}

\end{document}